\documentclass[11pt]{article}

\newcommand{\documentdate}{19 April 2017}

\usepackage{a4wide,latexsym,mystyle,graphicx,varioref,color,amsmath,cite}
\DeclareGraphicsExtensions{.eps,.jpg}
\usepackage{epstopdf}

\title{Partially separable convexly-constrained optimization \\
  with non-Lipschitzian singularities and its complexity}
\author{
  X. Chen\thanks{
    Department of Applied Mathematics, The Hong Kong Polytechnic University, Hong Kong.
    Email: \tt{xiaojun.chen@polyu.edu.hk}
  }
~ Ph. L. Toint\thanks{
    Namur Center for Complex Systems (naXys) and Department of Mathematics,
    University of Namur,
    61, rue de Bruxelles, B-5000 Namur, Belgium.
    Email: \tt{philippe.toint@unamur.be}}
~and H. Wang\thanks{
    Department of Applied Mathematics, The Hong Kong Polytechnic University, Hong Kong.
    Email: \tt{hong.wang@connect.polyu.hk}
  }
}
\date{\documentdate}
\newcommand{\ass}[2]{\label{ass-#1}
                     \begin{list}{}{\setlength{\leftmargin}{2.5cm}}
                     \item \hspace{-2.5cm} \framebox[2.0cm]{\bf #1} \,\,#2
                     \end{list}}

\DeclareMathOperator*\spanset{span}

\newcommand{\range}{{\rm range}}

\begin{document}

\comment{

\begin{titlepage}

\includegraphics[height=4cm]{/home/pht/Pictures/logos/UNamur.eps}
\hspace*{4cm}
\includegraphics[height=3cm]{/home/pht/Pictures/logos/naxys.eps}
\vspace*{3cm}
\begin{center}
\begin{minipage}[c]{12cm}
\vfill
\begin{center}
  {\sc
    Partially separable convexly-constrained optimization \\
   with non-Lipschitzian singularities and its complexity
}
\end{center}
\vfill
\centering{by X. Chen,  Ph. L. Toint and H. Wang}\\
\mbox{}
\vfill
\centering{Report NAXYS-??-2017 \hspace*{2 cm} \documentdate}\\*[1cm]
\vfill
\centering{\includegraphics[height=7cm]{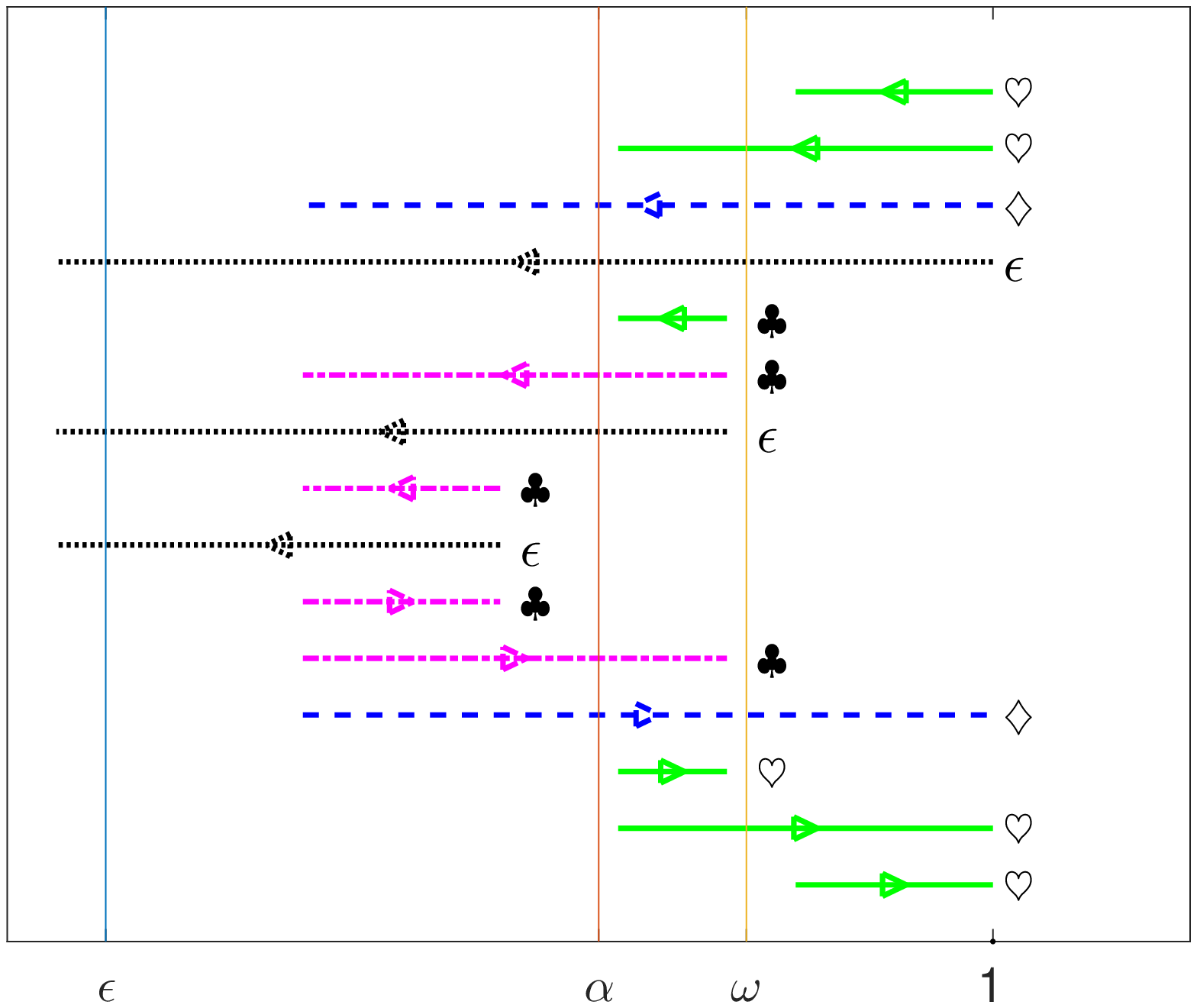}}
\end{minipage}
\end{center}
\vfill
\begin{center}
{\large
University of Namur, 61, rue de Bruxelles, B5000 Namur (Belgium)\\*[2ex]
{\tt http://www.unamur.be/sciences/naxys}}
\end{center}

\end{titlepage}
}

\maketitle

\begin{abstract}
\noindent
An adaptive regularization algorithm using high-order models is proposed for
partially separable convexly constrained nonlinear optimization problems
whose objective function contains non-Lipschitzian $\ell_q$-norm
regularization terms for $q\in (0,1)$.  It is shown that the algorithm using
an $p$-th order Taylor model for $p$ odd needs in general at most
$O(\epsilon^{-(p+1)/p})$ evaluations of the objective function and its
derivatives (at points where they are defined) to produce an
$\epsilon$-approximate first-order critical point. This result is obtained
either with Taylor models at the price of requiring the feasible set to
be 'kernel-centered' (which includes bound constraints and many other cases of
interest), or for non-Lipschitz models, at the price of passing the
difficulty to the computation of the step. Since this complexity bound
is identical in order to that already known for purely Lipschitzian
minimization subject to convex constraints \cite{CartGoulToin16a}, the new
result shows that introducing non-Lipschitzian singularities in the objective
function may not affect the worst-case evaluation complexity order. The
result also shows that using the problem's partially separable structure (if present)
does not affect complexity order either. A final (worse) complexity bound is derived
for the case where Taylor models are used with a general convex feasible set.
\end{abstract}

\vspace*{1cm}

\numsection{Introduction}

We consider the partially separable convexly constrained nonlinear optimization
problem:
\begin{equation}
\label{psprob}
  \min_{x \in \calF} f(x)
  = \sum_{i \in \calN} f_i(U_ix) + \sum_{i \in \calH} |U_ix|^q
  = \sum_{i \in \calN} f_i(x_i) + \sum_{i \in \calH} f_i(x_i)
\end{equation}
where $\calF$ is a non-empty closed convex set,
$\calN\cup\calH \eqdef \calM$,
$\calN\cap\calH = \emptyset$,
$f_i:\Re^n\rightarrow\Re$,
$q \in (0,1)$,
$f_i(x_i) = |x_i|^q= |U_ix|^q$ for $i \in \calH$ and where, for
$i \in \calM$, $x_i \eqdef U_ix$ with $U_i$ a (fixed) $n_i \times n$  matrix
with $n_i \leq n$.
Without loss of generality, we assume that $\|U_i\| = 1$ for all $i \in \calM$ and that
the ranges of the $U_i^T$ for $i\in \calN$ span $\Re^n$ in that the intersection
of the nullspaces of the $U_i$ is reduced to the origin\footnote{
  If the $\{U_i^T\}_{i \in \calN}$  do not span $\Re^n$, problem \req{psprob}
  can be modified without altering its optimal value by introducing an additional
  identically zero element term $f_0(U_0x)$ (say) in $\calN$ with associated
  $U_0$ such that $\cap_{i \in \calN}\ker(U_i) \subseteq {\rm
  range}(U_0^T)$.  It is clear that, since $f_0(x_0)=0$, it is
  differentiable with Lipschitz continuous derivative for any order $p\geq
  1$. Obviously, this covers the case where $\calN = \emptyset \neq \calH$.
}.
In what follows, the ``element functions'' $f_i$ ($i\in \calN$) will be
``well-behaved'' smooth functions with Lipschitz continuous
derivatives\footnote{Hence the symbol $\calN$ for ``nice''.}.
If $\calH \neq \emptyset$, we also require that
\beqn{UiH-conds}
n_i = 1
\eeqn
and (initially at least\footnote{We will drop this assumption in
Section~\ref{true-model-s}.}) that the feasible set is
'kernel centered', in the sense that, if $P_\calX[\cdot]$ is the orthogonal
projection ont the convex set $\calX$, then, for $i \in \calH$,
\beqn{Fbox}
P_{\ker(U_i)}[\calF] \subseteq \calF
\tim{ whenever } \ker(U_i) \cap \calF \neq \emptyset
\eeqn
in addition of $\calF$ being convex, closed and non-empty. As will be discussed
below (after Lemma~\ref{Lip-th}), we may assume without loss of generality
that, $\ker(U_i) \cap \calF \neq \emptyset$ (and thus $P_{\ker(U_i)}[\calF]
\subseteq \calF$) for all $i \in \calH$. 'Kernel centered' feasible sets
include boxes (corresponding to bound constrained problems), spheres/cylinders
centered at the origin and other sets such as
\beqn{kcex}
\left\{(x_1,x_2) \in \Re^{n_1+1} \mid  x_1 \in \calF_1 \tim{and}  g_1(x_1) \leq x_2  \leq -g_2(x_1)  \right\},
\ms
\calH = \{  2 \},
\eeqn
where $\calF_1$ is a non-empty closed convex set in $\Re^{n_1}$ and
$g_i(\cdot)$ are convex functions from $\Re^{n_1}$ to $\Re$  ($i=1,2$,
$n_1+1\leq n$) such that $g_i(x_1)\leq 0$ ($i=1,2$) for $x_1 \in \calF_1$.
Compositions using \req{kcex} recursively, rotations, cartesian products or
intersections of such sets are also kernel-centered.

Problem~\req{psprob} has many applications in engineering and science.  Using
the non-Lipschitz regularization function in the second term of the objective
function $f$ has remarkable advantages for the restoration of piecewise
constant images and sparse signals
\cite{BrucDonoElad09,BianChen15,LiuMaDaiZhan16}, and sparse variable
selection, for instance in bioinformatics
\cite{ChenNiuYuan13,HuanHoroMa18}. Theory and algorithms for solving
$q$-norm regularized optimization problems have been developed in
\cite{ChenGeYe14,ChenXuYe10,Lu14}.

The partially separable structure defined in problem \req{psprob} is
ubiquitous in applications of optimization. It is most useful in the frequent
case where $n_i \ll n$ and subsumes that of sparse optimization (in the
special case where the rows of each $U_i$ are selected rows of the identity
matrix). Moreover the decomposition in \req{psprob} has the advantage of being
invariant for linear changes of variables (only the $U_i$ matrices vary).

Partially separable optimization was proposed by Griewank and Toint
in \cite{GrieToin82a
}, studied for more than thirty years
(see \cite{
GoldWang93,
ChenDengZhan95,Gay96,ChenDengZhan98
,MareRichTaka14} for instance) and
extensively used in the popular {\sf CUTE}(st) testing environment
\cite{
GoulOrbaToin15b} as well as in the AMPL \cite{FourGayKern87}, {\sf LANCELOT}
\cite{ConnGoulToin92} and {\sf FILTRANE} \cite{GoulToin07b}
packages, amongst others.  In particular, the design of trust-region
algorithms exploiting the partially separable decomposition \req{psprob} was
investigated by Conn, Gould, Sartenaer and Toint in
\cite{ConnGoulSartToin96a,ConnGoulToin00} and Shahabuddin \cite{Shah96}.

Focussing now on the nice multivariate element functions ($i \in \calN$), we note that
using the partially separable nature of a function $f$ can be very useful if
one wishes to use derivatives of
\beqn{fN-def}
f_\calN(x) \eqdef \sum_{i \in \calN} f_i(U_ix) = \sum_{i \in \calN} f_i(x_i)
\eeqn
of order larger than one in the context of the $p$-th order Taylor series
\beqn{taylor}
T_{f_{\calN,p}}(x,s) = f_\calN(x) + \sum_{j=1}^p \frac{1}{j!}\nabla_x^jf_\calN(x)[s]^j.
\eeqn
Indeed, it may be verified that
\beqn{derjU}
\nabla_x^jf_\calN(x)[s]^j
= \sum_{i \in \calN} \nabla_{x_i}^jf_i(x_i)[U_is]^j.
\eeqn
This last expression indicates that only the $|\calN|$ tensors
$\{\nabla_{x_i}^jf_i(x_i)\}_{i\in \calN}$ of dimension $n_i^j$ needs to be computed
and stored, a very substantial gain compared to the $n^j$-dimensional
$\nabla_x^jf_\calN(x)$ when (as is common) $n_i \ll n$ for all $i$.  It may
therefore be argued that exploiting derivative tensors of order larger than 2
--- and thus using the high-order Taylor series \req{taylor} as a local model
of $f(x+s)$ in the neighbourhood of $x$ --- may be practically feasible if $f$ is
partially separable. Of course the same comment applies to
\beqn{fH-def}
f_\calH(x) \eqdef \sum_{i \in \calH} f_i(U_ix) = \sum_{i \in \calH} f_i(x_i)
\eeqn
whenever the required derivatives of $f_i(x_i) = |x_i|^p$ ($i\in \calH$) exist.

Interestingly, the use of high-order Taylor models for optimization was
recently investigated by Birgin \emph{et al.} \cite{BirgGardMartSantToin17}
in the context of adaptive regularization algorithms for unconstrained
problems.  Their proposal belongs to this emerging class of methods pioneered
by Griewank \cite{Grie81}, Nesterov and Polyak \cite{NestPoly06} and Cartis,
Gould and Toint \cite{CartGoulToin11,CartGoulToin11d} for the unconstrained
case and by these last authors in \cite{CartGoulToin12b} for the convexly
constrained case of interest here.  Such methods are distinguished by their
excellent evaluation complexity, in that they need at most
$O(\epsilon^{-(p+1)/p})$ evaluations of the objective function and their
derivatives to produce an $\epsilon$-approximate first-order critical point,
compared to the $O(\epsilon^{-2})$ evaluations which might be necessary for
the steepest descent and Newton's methods (see \cite{Nest04} and
\cite{CartGoulToin10a} for details).  However, most adaptive regularization
methods rely on a non-separable regularization term in the model of the
objective function, making exploitation of structure difficult\footnote{The
only exception we are aware of is the unpublished note
\cite{GoulHoggReesScot16} in which a $p$-th order Taylor model is coupled
with a regularization term involving the (totally separable) $q$-th power of
the $q$ norm ($q \geq 1$).}.

The purpose of the present paper is twofold.  Its first
aim is to show that worst-case evaluation complexity for nonconvex
minimization subject to convex constraints is not affected by the introduction
of  non-Lipschitzian singularities in the objective function.
The second and concurrent one is to show that this complexity is not affected
either by the use of partially separable structure, if present in the problem.

The remaining of the paper is organized as follows. Section~\ref{optimality-s}
establishes a necessary first-order optimality condition for the
non-Lipschitzian case. Section~\ref{psarp-s} then
introduces the partially separable adaptive regularization algorithm for this
problem while Section~\ref{complexity-s} is devoted to its worst-case evaluation
complexity analysis for the case where Taylor models are used with a
kernel-centered feasible set.  Section~\ref{true-model-s} drops the
kernel-centered assumption for non-Lipschitz models and Taylor models.
The results are discussed in Section~\ref{discuss-s} and
some final conclusions and perspectives are presented in Section~\ref{concl-s}.

\noindent
{\bf Notations.}
In what follow, $\|x\|$ denotes the Euclidean norm of the vector $x$ and $\|T\|_p$
the recursively induced Euclidean norm on the $p$-th order tensor $T$ (see
\cite{BirgGardMartSantToin17,CartGoulToin16a} for details).  The notation
$T[s]^i$ means that the tensor $T$ is applied to $i$ copies of
the vector $s$. For any set $\calX$, $|\calX|$ denotes its
cardinality. For any $\calI \subseteq \calM$, we also denote
$f_\calI(x)= \sum_{i\in \calI}f_i(x)$.

\numsection{First-order necessary conditions}\label{optimality-s}

In this section, we first present exact and approximate first-order necessary
conditions for a local  minimizer of problem \req{psprob}.
Such conditions for optimization problems with non-Lipschitzian
singularities have been independently defined in the scaled form
\cite{ChenXuYe10} or in subspaces \cite{BianChen17,ChenNiuYuan13}. In a recent paper
\cite{ChenLeiLuYe17}, KKT necessary optimality conditions for constrained
optimization problems with non-Lipschitzan singularities are studied under the
relaxed constant positive linear dependence and basic qualification.
The above optimality conditions take the singularity into account
by no longer requiring that the gradient (for unconstrained problems, say) nearly
vanishes at an approximate solution $x_\epsilon$ (which would be impossible if the
singularity is active) but by requiring that a scaled version of this
requirement holds in that $\|X_\epsilon\nabla_x^1f(x_\epsilon)\|$ is suitably
small, where $X_\epsilon $ is a diagonal matrix whose diagonal entries are the
components of $x_\epsilon$. Unfortunately, if the $i$-th component of
$x_\epsilon$ is small but not quite small enough to consider that the
singularity is active for variable $i$ (say it is equal to $2\epsilon$), the
$i$-th component of $\nabla_x^1f(x)$ can be as large as a multiple
of $\epsilon^{-1}$. As a result, comparing worst-case evaluation complexity
bounds with those known for purely Lipschitz continuous problems (such as
those proposed in \cite{BirgGardMartSantToin17} or \cite{CartGoulToin16a})
may be misleading, since these latter conditions would never accept an
approximate first-order critical point with such a large gradient. In order to
avoid these pitfalls, we now propose a stronger definition of approximate
first-order critical point for non-Lipschitzian problems where such
``border-line'' situations do not occur. The new definition is also makes use of
subspaces but exactly reduces to the standard condition for Lipschitzian
problems if the singularity is not active at $x_\epsilon$, even if it is close
to it.

Given a vector $x\in \Re^n$ and $\epsilon\ge0$, denote
\begin{eqnarray*}\label{CR-defs}
  \calC(x,\epsilon)\eqdef \{i\in\calH\ |\ |U_i x| \le \epsilon\}, \quad
  \calR(x,\epsilon)\eqdef \bigcap_{i\in\calC(x,\epsilon)}\ker(U_i)=\left[
  \spanset_{i\in\calC{(x,\epsilon)}}\{U_i^T\}\right]^{\bot}
\end{eqnarray*}
and
\begin{equation*}
  \calW(x,\epsilon)\eqdef \calN \cup (\calH\setminus \calC(x,\epsilon)).
\end{equation*}

For convenience, if $\epsilon=0$, we denote $\calC(x)\eqdef \calC(x,0)$, $\calR(x)
\eqdef \calR(x,\epsilon)$ and $\calW(x)\eqdef \calW(x,0)$.

Observe that the definition of $\calR(x,\epsilon)$ above gives that
\beqn{RperpinspanH}
\calR(x,\epsilon)^\perp \subseteq \spanset_{i \in \calH}\{U_i^T\}.
\eeqn
Also note that any $x\in\Re^n$ can be decomposed uniquely as $x=y+z$ where
$y\in\calR(x)^{\bot}$ and $z\in\calR(x)$. By the definition of $\calR(x)$, it
is not difficult to verify that
\[
U_i z=0,~\forall i\in \calC(x)
\quad\mbox{and}\quad
x \in\calR(x).
\]
Finaly note that, although $f(x)$ is nonsmooth if $\calH \neq \emptyset$,
$f_{\calW(x,\epsilon)}(x)$ is as differentiable as the $f_i(x)$ for
$i\in\calN$ and any $\epsilon \geq 0$.  This allows us to formulate our
first-order necessary condition.

\lthm{lemma:first-order-necessary}{
If $x_*\in \calF$ is a local minimizer
of problem \req{psprob}, then
\begin{equation}\label{eq:first-order-condition}
  \chi_f(x_*) =0,
\end{equation}
where, for any $x\in \calF$,
\beqn{chi-def}
\chi_f(x_*)= \chi_f(x_*,0)
\eqdef
  \left|\min_{\mystack{x+d\in\calF}{d\in\calR(x),\|d\|\leq 1} }
      \nabla_x^1 f_{\calW(x)}(x)^T d\right|.
\eeqn
}

\proof{
Suppose first that $\calR(x_*) = \{0\}$ (which happens if $x_* = 0 \in \calF$ and
$\spanset_{i\in \calH}\{U_i^T\} = \Re^n$).  Then \req{eq:first-order-condition}-\req{chi-def}
holds vacuously. Now suppose that
$\calR(x_*)$ contains at least one nonzero element.  By assumption, there
exists $\delta_{x_*}>0$ such that
\begin{equation*}
\begin{split}
f(x_*)&=\min\{f_{\calN}(x)+f_{\calH}(x)\ |\ x\in {\cal F}, ~\|x-x_*\|\leq\delta_{x_*}\} \\
&=\min\{f_{\calN}(y+z)+f_{\calH}(y+z)\ |\  y+z\in {\cal F},y\in\calR(x_*)^{\bot},z\in\calR(x_*),
 \|y+z-x_*\|\le \delta_{x_*}\}\\
&\le\min\{f_{\calN}(y+z)+\sum_{i\in\calH}|U_i(y+z)|^q\ |\
 y+z\in {\cal F},y=0,z\in\calR(x_*),\|z-x_*\|\le\delta_{x_*}\} \\
 &=\min\{f_{\calN}(z)+\sum_{i\in\calH}|U_i z|^q\ |\ z\in {\cal F}\cap\calR(x_*),
 \|z-x_*\|\leq \delta_{x_*}\}\\
 &=\min\{f_{\calN}(z)+\sum_{i\in\calH\backslash\calC({x_*})}|U_i z|^q\ |\
 z\in {\cal F}\cap\calR(x_*),\|z-x_*\|\leq \delta_{x_*}\}.
\end{split}
\end{equation*}

We now introduce a new problem, which is problem \req{psprob} reduced to
$\calR(x_*)$, namely,
\begin{equation}\label{eq:reduced-problem-of-original-prob}
  \begin{split}
    \min&\quad f_{\calW(x_*)}(z)=f_{\calN}(z)+\sum_{i\in \calH\backslash\calC(x_*)}
    |U_i z|^q,\\
  \mbox{s.t.} &\quad z\in {\cal F}\cap\calR(x_*)
  \end{split}
\end{equation}
where the gradient $\nabla_z^1 f_{\calW(x_*)}(z)$ is locally Lipschitz continuous
in some (bounded) neighborhood of $x_*$.
It then follows from $x_*\in\calR(x_*)$ that
\begin{equation*}
f_{\calW(x_*)}(x_*)=f_{\calN}(x_*)+\sum_{i\in\calH\backslash\calC(x_*)}|U_i x_*|^q
          =f(x_*).
\end{equation*}
Therefore, we have that
\begin{eqnarray*}
  f_{\calW(x_*)}(x_*)
  \leq \min\{f_{\calW(x_*)}(z)\ |\ z\in {\cal F}\cap \calC(x_*),\|z-x_*\|\leq\delta_{x_*}\}
\end{eqnarray*}
which implies that $x_*$ is a local minimizer of problem
\req{eq:reduced-problem-of-original-prob}. Hence, we have
\begin{equation}\label{first-order-v}
\nabla_z^1 f_{\calW(x_*)}(x_*)^T(z-x_*)\geq 0, \quad z\in {\cal F}\cap \calR(x_*).
\end{equation}
In addition,
\[
\{d = 0\}\subseteq\{d\ \mid\ x_*+d \in \calF,d\in\calR(x_*),\|d\|\le1\}
\subseteq\{d\ \mid\ x_*+d \in \calF, d\in \calR(x_*)\}
\]
which gives the desired result \req{eq:first-order-condition}-\req{chi-def}.
} 

\noindent
We call $x_*$ a \emph{first-order stationary point} of \req{psprob}, if $x_*$
satisfies the relation \req{eq:first-order-condition}  in
Theorem~\ref{lemma:first-order-necessary}.
For $\epsilon>0$, we call $x_\epsilon$ an \emph{$\epsilon$-approximate
  first-order stationary point} of \req{psprob}, if $x_\epsilon$ satisfies
\begin{equation}\label{eq:epsilon-approximate-first-order-stationary}
  \chi_f(x_\epsilon,\epsilon)
  \eqdef \left|\min_{\mystack{x+d\in\calF}{d \in\calR(x_\epsilon,\epsilon),\|d\|\leq 1}}
          \nabla_x^1f_{\calW(x_\epsilon,\epsilon)}(x_\epsilon)^T d\right|
  \le\epsilon.
\end{equation}

\lthm{prop:epsilon-scaled-and-epsilon-stationary}{
  Let $x_\epsilon$ be an $\epsilon$-approximate first-order stationary point
  of \req{psprob}. Then any cluster point of $\{x_\epsilon\}_{\epsilon>0}$ is
  a first-order stationary point of problem \req{psprob} as $\epsilon\rightarrow0$.
}

\proof{
Suppose that $x_*$ is any cluster  point of $\{x_\epsilon\}_{\epsilon>0}$.
Hence there must exist an infinite sequence $\{\epsilon_k\}$ converging to
zero and an infinite sequence $\{x_{\epsilon_k}\}_{k\geq 0}\subseteq\{x_\epsilon\}_{\epsilon>0}$
such that $x_* = \lim_{k\rightarrow\infty}x_{\epsilon_k}$ and
$x_{\epsilon_k}$ is an $\epsilon_k$-approximate
first-order stationary point of \req{psprob} for eack $k\geq 0$.
If $\calR(x_*)=\{0\}$, \req{eq:first-order-condition} holds vacuously and
hence $x_*$ is a first-order stationary point. Suppose therefore that
$\calR(x_*)$ contains at least one nonzero element, implying that the dimension of
$\calR(x_*)$ is strictly positive.

First of all, we claim that there must exist $k_*\geq 0$
such that $\calR(x_{\epsilon_k},\epsilon_k)^{\bot}\subseteq\calR(x_*)^{\bot}$
for any $k\ge k_*$. Indeed, if that is not the case, there exists a
subsequence of $\{x_{\epsilon_k}\}$, say $\{x_{\epsilon_{k_j}}\}$, such that
$\lim_{j\rightarrow\infty}\epsilon_{k_j}=0$ and
$\calR(x_{\epsilon_{k_j}},\epsilon_{k_j})^{\bot}\not\subseteq\calR(x_*)^{\bot}$
for all $j$. Using now \req{RperpinspanH} and the fact that $\calH$ is a
finite set, we obtain that there must exist an $i_0\in\calH$ such
that $i_0\in\calC(x_{\epsilon_{k_{j_t}}},\epsilon_{k_{j_t}})$ but $i_0\notin
\calC(x_*)$ where $\{k_{j_t}\}\subseteq\{k_j\}$ with $t=1,2,\cdots$. For
convenience, we continue to use $\{k_j\}$ to denote its subsequence
$\{k_{j_t}\}$. Hence, we have that
\begin{eqnarray*}
|U_{i_0} x_{\epsilon_{k_j}}|\leq \epsilon_{k_j}.
\end{eqnarray*}
Let $j$ go to infinity. It then follows from the above inequality that
$|U_{i_0} x_*|=0$, which contradicts the fact that $i_0\notin\calC(x_*)$.
Thus, we conclude that, for some $k_*\geq 0$ and all $k\ge k_*$,
$\calR(x_{\epsilon_k},\epsilon_k)^{\bot}\subseteq\calR(x_*)^{\bot}$.
Therefore we have that
$\calR(x_*)\subseteq\calR(x_{\epsilon_k},\epsilon_k)$ for $k \geq k_*$.

For any fixed $\epsilon_k$ approximate first-order stationary point $x_{\epsilon_k}$,
consider the following two minimization problems.
\begin{equation}\label{eq:two-constrained-linear-minimization-prob-a}
\begin{split}
  \min& \quad \nabla_x^1 f_{\calW(x_{\epsilon_k},\epsilon_k)}(x_{\epsilon_k})^T d,  \\
  \mbox{s.t.} & \quad x_{\epsilon_k}+d\in\calF,d\in\calR(x_{\epsilon_k},\epsilon_k), \|d\|\le 1,
\end{split}
\end{equation}
and
\begin{equation}\label{eq:two-constrained-linear-minimization-prob-b}
 \begin{split}
 \min & \quad \nabla_x^1f_{\calW(x_{\epsilon_k},\epsilon_k)}(x_{\epsilon_k})^T d, \\
 \mbox{s.t.}  & \quad x_{\epsilon_k}+d\in\calF,d\in\calR(x_*), \|d\|\leq 1.
 \end{split}
\end{equation}

Since $d=0$ is a feasible point of both problems
\req{eq:two-constrained-linear-minimization-prob-a} and
\req{eq:two-constrained-linear-minimization-prob-b}, the minimum values of
\req{eq:two-constrained-linear-minimization-prob-a} and
\req{eq:two-constrained-linear-minimization-prob-b} are both nonpositive.
Moreover, it follows from $\calR(x_*)\subseteq\calR(x_{\epsilon_k},\epsilon_k)$
that the minimum value of \req{eq:two-constrained-linear-minimization-prob-b}
is not smaller than that of \req{eq:two-constrained-linear-minimization-prob-a}.

Hence, from \req{eq:epsilon-approximate-first-order-stationary}, we have that for any
$x_{\epsilon_k}$,
\begin{equation}\label{eq:epsilon-first-order-point-k-defn}
  \left|\min_{\mystack{x_{\epsilon_k}+d\in\calF}{d\in\calR(x_*), \|d\|\leq 1}}
  \nabla_x^1f_{\calW(x_{\epsilon_k},\epsilon_k)}(x_{\epsilon_k})^T d\right|
 \le
 \left|\min_{\mystack{x_{\epsilon_k}+d\in\calF}{d\in\calR(x_{\epsilon_k},\epsilon_k), \|d\|\leq 1}}
     \nabla_x^1f_{\calW(x_{\epsilon_k},\epsilon_k)}(x_{\epsilon_k})^T d\right|\le \epsilon_k.
\end{equation}
Suppose that $d_{\epsilon_k}$ is a  minimizer of problem
\req{eq:two-constrained-linear-minimization-prob-b}, then
\req{eq:epsilon-first-order-point-k-defn} implies that
\begin{equation}\label{eq:epsilon-first-order-point-k-defn-imply}
  -\epsilon_k \leq \nabla_x^1f_{\calW(x_{\epsilon_k},\epsilon_k)}(x_{\epsilon_k})^T d_{\epsilon_k} \leq 0,
\end{equation}
where $d_{\epsilon_k}$ should satisfy that $x_{\epsilon_k}+d_{\epsilon_k}\in\calF$,
$d_{\epsilon_k}\in \calR(x_*)$ and $\|d_{\epsilon_k}\|\leq 1$.
Note that, since $d_{\epsilon_k}\in \calR(x_*)$,
\begin{equation}\label{eq:reduce-dk-subspace-hk}
\begin{split}
\nabla_x^1f_{\calW(x_{\epsilon_k},\epsilon_k)}(x_{\epsilon_k})^T d_{\epsilon_k} & = \left(\nabla_x f_{\calN}(x_{\epsilon_k})
+\sum_{i\in\calH\setminus\calC(x_{\epsilon_{k}})}q|U_i x_{\epsilon_k}|^{q-1}\sign(U_i x_{\epsilon_k})U_i^T \right)^T d_{\epsilon_k}  \\
& = \nabla_x f_{\calN}(x_{\epsilon_k})^T d_{\epsilon_k}+\sum_{i\in\calH\setminus\calC(x_{\epsilon_k})}q|U_i x_{\epsilon_k}|^{q-1}\sign(U_i x_{\epsilon_k})U_i d_{\epsilon_k} \\
& = \nabla_x f_{\calN}(x_{\epsilon_k})^T d_{\epsilon_k}+\sum_{i\in\calH}q|U_i x_{\epsilon_k}|^{q-1}\sign(U_i x_{\epsilon_k})U_i d_{\epsilon_k}.
\end{split}
\end{equation}

From the compactness of $\{d\ |\ \|d\|\le1\}$, we know that there must exist a subsequence
of $\{d_{\epsilon_k}\}$ such that $d_{\epsilon_{k_j}}\rightarrow d_*\in\calR(x_*)$
with $\|d_*\|\leq 1$ as $j$ goes to infinity.
Since for $i\in\calH\backslash\calC(x_*)$, we have
$\lim_{k\rightarrow\infty}|U_i x_{\epsilon_k}|^{q-1}=|U_ix_*|^{q-1}$.
Let $k$ go to infinity in \req{eq:epsilon-first-order-point-k-defn-imply}
and \req{eq:reduce-dk-subspace-hk}, and we obtain that
\begin{equation*}
  0=\nabla_x^1f_{\calW(x_{\epsilon_k},\epsilon_k)}(x_*)^T d_* = \nabla f_{\calN}(x_*)^T d_*
  + \sum_{i\in\calH\backslash\calC(x_*)} q|U_i x_*|^{q-1}\sign(U_ix_*)U_id_*,
\end{equation*}
which implies that
\[
\min_{\mystack{x_*+d\in\calF}{d\in\calR(x_*), \|d\|\leq 1}} \nabla_x^1f_{\calW(x_{\epsilon_k},\epsilon_k)}(x_*)^T d=\nabla_x^1f_{\calW(x_*)}(x_*)^T d_*=0
\]
and completes the proof.
}

\numsection{A partially separable regularization algorithm}\label{psarp-s}

We now examine the desired properties of the element functions $f_i$
more closely.
Assume first that, for $i \in \calN$, each element
function $f_i$ is $p$ times continuously differentiable
and its $p$-th order derivative tensor $\nabla_x^p f_i$ is globally
Lipschitz continuous with constant $L_i\geq 0$ in the sense that, for all
$x_i,y_i \in \range(U_i)$,
\beqn{tensor-Lip-fi}
\|\nabla_{x_i}^pf_i(x_i) - \nabla_{x_i}^p f_i(y_i)\|_p \leq L_i \|x_i-y_i\|.
\eeqn
It can be shown (see \req{fi-Lip} below) that this assumption implies that,
for $i \in \calN$,
\beqn{f-Lip-1}
f_i(x_i+s_i) = T_{f_i,p}(x_i,s_i) + \frac{1}{(p+1)!} \tau_i L_i \|s_i\|^{p+1}
\tim{ with }| \tau_i | \leq 1,
\eeqn
where $s_i = U_is$.

Because the quantity $\tau_i L_i$ in \req{f-Lip-1} is usually unknown in
practice, it is impossible to use \req{f-Lip-1} directly to model the objective
function in a neighbourhood of $x$.  However, we may replace this term with an
adaptive parameter $\sigma_i$, which yields the following
$(p+1)$-th order model for the $i$-th element ($i \in \calN$):
\beqn{model-fi}
m_i(x_i,s_i)=T_{f_i,p}(x_i,s_i)+\frac{1}{(p+1)!}\ \sigma_i\|s_i\|^{p+1}.
\eeqn

There is more than one possible choice for defining the element models for $i
\in \calH$. The first\footnote{Another choice is discussed in
  Section~\ref{true-model-s}.}
is to pursue the line of polynomial Taylor-based models,
for which we need the following technical result.

\llem{combined-model-l}{
  We have that, for $i \in \calH$ and all $x,s \in \Re^n$ with $U_ix \neq 0
  \neq U_i(x+s)$,
  \beqn{combined-model}
  |x_i+s_i|^q
  = |x_i|^q + q \sum_{j=1}^\infty \frac{1}{j!}\left(\prod_{\ell=1}^{j-1}(q-\ell)\right)
    |x_i|^{q-j} \mu(x_i,s_i)^j,
  \eeqn
  where
  \beqn{mu-def}
  \mu(x_i,s_i) \eqdef \left\{ \begin{array}{ll}
    ~~\,s_i      & \tim{if} x_i > 0 \tim{and} x_i+s_i > 0, \\
   -s_i          & \tim{if} x_i < 0 \tim{and} x_i+s_i < 0, \\
    -(2x_i+s_i) & \tim{if} x_i > 0 \tim{and} x_i+s_i < 0, \\
   ~~\,2x_i+s_i & \tim{if} x_i < 0 \tim{and} x_i+s_i > 0.
  \end{array} \right.
  \eeqn
}

\proof{
  If $y \in \Re_+$, it can be verified that the Taylor expansion $|y+z|^q$ at
  $y\neq 0$ and $y+z \in \Re_+$ is given by
  \beqn{tayl-y}
  [y+z]^q
  = y^q + q\sum_{j=1}^\infty \frac{1}{j!}  \left[\prod_{\ell=1}^{j-1}(q-\ell)\right] y^{q-j} z^j.
  \eeqn
  Let us now consider $i\in \calH$.
  Relation \req{tayl-y} yields that, if $x_i>0$ and $x_i+s_i > 0$,
  \beqn{mod++}
  |x_i+s_i|^q
  = |x_i|^q + q \sum_{j=1}^\infty \frac{1}{j!}\left[\prod_{\ell=1}^{j-1}(q-\ell)\right]
   |x_i|^{q-j} s_i^j.
  \eeqn
  By symmetry, if we have that if $x_i <0$ and $x_i+s_i <  0$, then
  \beqn{mod--}
  |x_i+s_i|^q
  = |x_i|^q + q \sum_{j=1}^\infty \frac{(-1)^j}{j!}\left[\prod_{\ell=1}^{j-1}(q-\ell)\right]
   |x_i|^{q-j} s_i^j.
   \eeqn
  Moreover, if $x_i >0$ and $x_i+s_i < 0$, then
  \beqn{mod+-}
  |x_i+s_i|^q
  = |-x_i|^q + q \sum_{j=1}^\infty \frac{(-1)^j}{j!}\left[\prod_{\ell=1}^{j-1}(q-\ell)\right]
   |-x_i|^{q-j} (2x_i+s_i)^j.
  \eeqn
  Symmetrically, if $x_i <0$ and $x_i+s_i > 0$, then
  again,
  \beqn{mod-+}
  |x_i+s_i|^q
  = |-x_i|^q + q \sum_{j=1}^\infty \frac{1}{j!}\left[\prod_{\ell=1}^{j-1}(q-\ell)\right]
  |-x_i|^{q-j} (2x_i+s_i)^j
   \eeqn
   \req{combined-model}-\req{mu-def} then trivially results from
   \req{mod++}-\req{mod-+} and the identity $|-x_i| = |x_i|$.
} 

\begin{figure}[htbp]
\begin{center}
\vspace*{1.5mm}
\includegraphics[height=6.5cm]{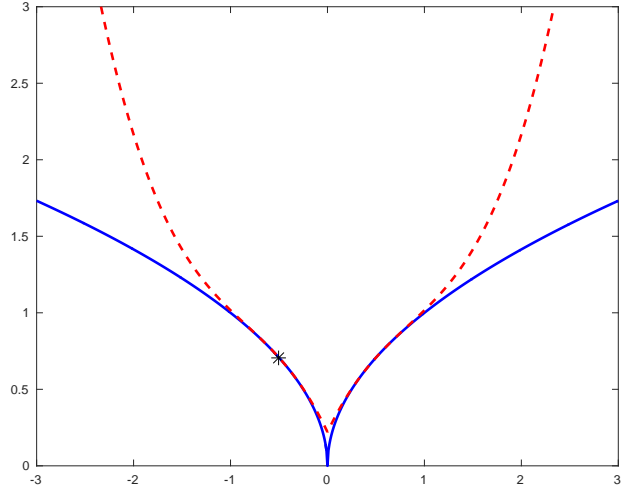} 
\caption{\label{twosided-fig}The square root function (continuous) and its
  two-sided model with $p=3$ evaluated at $x_i=-\half$ (dashed)}
\end{center}
\vspace*{-5mm}
\end{figure}

\noindent
We now slightly abuse notation by defining
\beqn{modiH}
T_{|\cdot|^q,p}(x_i,s_i)
\eqdef \left\{\begin{array}{ll}
T_{x^q,p}(x_i, s_i)
& \tim{if }  x_i >0 \tim{and} x_i+s_i > 0, \\*[2ex]
T_{(-x)^q,p}(x_i, -s_i)
& \tim{if }  x_i <0 \tim{and} x_i+s_i < 0, \\*[2ex]
T_{(-x)^q,p}(-x_i, 2x_i+s_i)
& \tim{if } x_i > 0 \tim{and} x_i+s_i < 0, \\*[2ex]
T_{x^q,p}(-x_i, 2x_i+s_i)
& \tim{if } x_i <0 \tim{and} x_i+s_i>0.
\end{array}\right.
\eeqn
We are now in position to define the regularized ``two-sided'' model for the
element function $f_i$ ($i\in \calH$) as
\beqn{mod-hir}
m_i(x_i,s_i)
\eqdef T_{|\cdot|^q,p}(x_i,s_i).
\eeqn
Figure~\ref{twosided-fig} illustrates the two-sided model
\req{modiH}-\req{mod-hir} for $x_i= - \half$, $p=3$, $q=\half$.

We may now build the complete  model for $f$ at $x$ as
\beqn{full-model}
m(x,s) =  \bigsum_{i\in\calM}  m_i(x_i,s_i).
\eeqn

The algorithm considered in this paper exploits the model \req{full-model} as
follows.  At each iteration $k$, the model \req{full-model} taken at the
iterate $x=x_k$ is (approximately) minimized in order to define a step
$s_k$. If the decrease in the objective function value along $s_k$ is
comparable to that predicted by the Taylor model, the trial point $x_k+s_k$ is
accepted as the new iterate and the regularization parameters $\sigma_{i,k}$
(i.e. $\sigma_i$ at iteration $k$) possibly updated. The
process is terminated when an approximate local minimizer is found, that is
when, for some $k \geq 0$,
\beqn{chi-eps}
\chi_f(x_k,\epsilon) \leq \epsilon.
\eeqn

In order to simplify notation in what follows, we make the following
definitions:
\[
\calC_k \eqdef \calC(x_k,\epsilon),
\ms
\calR_k \eqdef \calR(x_k,\epsilon),
\ms
\calW_k \eqdef \calW(x_k,\epsilon),
\]
and
\[
\calC_k^+ \eqdef \calC(x_k+s_k,\epsilon),
\ms
\calR_k^+ \eqdef \calR(x_k+s_k,\epsilon),
\ms
\calW_k^+ \eqdef \calW(x_k+s_k,\epsilon).
\]

Having defined the criticality measure \req{chi-def}, it is natural to use
this measure also for terminating the approximate
model minimization: to find $s_k$, we therefore minimize $m(x_k,s)$ over $s$
until, for some constant $\theta \geq 0$ and some exponent $r > 1$,
\beqn{mod-term}
\chi_m(x_k,s_k,\epsilon)=\chi_{m_{\calW_k^+}}(x_k,s_k,\epsilon)
\leq \min\left[ \, \quarter q^2 \min_{i\in \calH \cap \calW_k^+}|U_i(x_k+s_k)|^r,
             \, \theta \|s_k\|^p \, \right]
\eeqn
where
\beqn{chim-def}
\chi_{m_{\calW_k^+}}(x_k,s_k,\epsilon)
\eqdef \left| \min_{\mystack{x_k+s_k+d \in \calF}{ d \in \calR_k^+, \|d\|\leq 1}}
                                    \nabla_s^1m_{\calW_k^+}(x_k,s_k)^Td \right|.
\eeqn

We also require that, once $|U_i(x_k+s)| < \epsilon$ occurs for some
$i \in \calH$ in the course of the model minimization, it is fixed at
this value, meaning that the remaining minimization is carried out in
$\calR(x_k+s, \epsilon)$.  Thus the dimension of $\calR(x_k+s,\epsilon)$ (and
therefore of $\calR(x_k,\epsilon)$) is monotonically non-increasing during the
step computation and across iterations.
Note that computing a step $s_k$ satisfying \req{mod-term}
is always possible since the subspace $\calR(x_k+s,\epsilon)$ can only become
smaller during the model minimization and since we have seen in
Section~\ref{optimality-s} that $\chi_m(x_k,s_k)=0$ at any local
minimizer of $m_{\calW(x_k+s,\epsilon)}(x_k,s)$.

\subsection{The algorithm}

We now introduce some notation useful for describing our algorithm.
Define
\[
x_{i,k} \eqdef U_ix_k,
\ms
s_{i,k} \eqdef U_is_k.
\]
Also let
\[
\delta f_{i,k} \eqdef f_i(x_{i,k}) - f_i(x_{i,k}+s_{i,k})
\]
\[
\delta f_k
\eqdef f_{\calW_k^+}(x_k) - f_{\calW_k^+}(x_k+s_k)
= \sum_{i \in \calW_k^+} \delta f_{i,k},
\]
\[
\delta m_{i,k} \eqdef m_i(x_{i,k},0) - m_i(x_{i,k},s_{i,k}),
\]
\[
\delta m_k
\eqdef m_{\calW_k^+}(x_k,0) - m_{\calW_k^+}(x_k,s_k)
= \sum_{i \in \calW_k^+} \delta m_{i,k},
\]
and
\beqn{deltaT-def}
\begin{array}{lcl}
\delta T_k
& \eqdef &  T_{f_{\calW_k^+},p}(x_k,0) - T_{f_{\calW_k^+},p}(x_k,s_k) \\*[2ex]
&    =   & [T_{f_{\calN},p}(x_k,0) - T_{f_{\calN},p}(x_k,s_k)]
+ [T_{|\cdot|_{\calH\setminus\calC_k^+,p}}(x_k,0)
  - T_{|\cdot|_{\calH\setminus\calC_k^+,p}}(x_k,s_k)]\\*[2ex]
&    =   & \delta m_k+\bigfrac{1}{(p+1)!}\bigsum_{i \in \calN}\sigma_{i,k}\|s_{i,k}\|^{p+1}.
\end{array}
\eeqn
The partially separable adaptive regularization algorithm is now formally
stated as Algorithm~\vref{psarp-algo}.

\algo{psarp-algo}{Partially Separable Adaptive Regularization}{
\begin{description}
\item[Step 0: Initialization:] $x_0\in\calF$ and $\{\sigma_{0,i}\}_{i\in\calN}
  >0$ are given as well as the accuracy $\epsilon \in (0,1]$ and constants $0
<\gamma_0 < 1 <   \gamma_1 \leq \gamma_2$, $\eta \in (0,1)$, $\theta \geq0$,
$\sigma_{\min} \in  (0,\min_{i\in\calN}\sigma_{0,i}]$ and  $\kap{big}>1$. Set $k =0$.
\item[Step 1: Termination:] Evaluate $f(x_k)$ and
  $\{\nabla_x^1 f_{\calW_k}(x_k)\}$. If
  $\chi_f(x_k,\epsilon) \leq \epsilon$, return $x_\epsilon=x_k$ and
  terminate. Otherwise evaluate
  $\{\nabla_x^i f_{\calW_k}(x_k)\}_{i=2}^p$.
\item[Step 2: Step computation:] Compute a
  step $s_k\in \calR_k$ such that $x_k+s_k \in \calF$, $m(x_k,s_k)<  m(x_k,0)$
  and \req{mod-term} holds.
\item[Step 3: Step acceptance:] Compute
\beqn{rhok-def}
\rho_k = \bigfrac{\delta f_k}{\delta T_k}
\eeqn
and set $x_{k+1} = x_k$ if $\rho_k < \eta$, or $x_{k+1} = x_k+s_k$ if $\rho_k \geq \eta$.

\item[Step 4: Update the ``nice'' regularization parameters:] For
$i \in \calN$, if
\beqn{up-cond}
f_i(x_{i,k}+s_{i,k}) > m_i(x_{i,k},s_{i,k})
\eeqn
set
\beqn{sig-incr}
\sigma_{i,k+1} \in [ \gamma_1 \sigma_{i,k}, \gamma_2 \sigma_{i,k} ].
\eeqn
Otherwise, if either
\beqn{muchover-neg}
\rho_k \geq \eta
\tim{ and }
\delta f_{i,k} \leq 0
\tim{ and }
\delta f_{i,k} < \delta m_{i,k} - \kap{big}|\delta f_k|
\eeqn
or
\beqn{muchover-pos}
\rho_k \geq \eta
\tim{ and }
\delta f_{i,k} > 0
\tim{ and }
\delta f_{i,k} > \delta m_{i,k} + \kap{big}|\delta f_k|
\eeqn
then set
\beqn{sig-decr}
\sigma_{i,k+1} \in [\max[\sigma_{\min}, \gamma_0\sigma_{i,k}], \sigma_{i,k} ],
\eeqn
else set
\beqn{sig-keep}
\sigma_{i,k+1} = \sigma_{i,k}.
\eeqn
Increment $k$ by one and go to Step~1.
\end{description}
}

Note that an $x_0 \in \calF$ can always be computed by projecting an
infeasible starting point onto $\calF$.
The idea of the second and third parts of \req{muchover-neg} and
\req{muchover-pos} is to identify cases where the model $m_i$
overestimates the element function $f_i$ to an excessive extent, leaving some
space for reducing the regularization and hence allowing longer steps.
The requirement that $\rho_k \geq \eta$ in both \req{muchover-neg} and
\req{muchover-pos} is intended to prevent a situation where a particular
regularization parameter is increased and another decreased at a given
unsuccessful iteration, followed by the opposite situation at the next
iteration, potentially leading to cycling. Other more elaborate mechanisms can
be designed to achieve the same goal, such as attempting to reduce a given
regularization parameter at most a fixed number of times before the occurence
of a successful iteration, but we do not investigate those alternatives in
detail here. The idea of the second and third parts of \req{muchover-neg} and
\req{muchover-pos} is simply to identify cases where the model $m_i$
overestimates the element function $f_i$ to an excessive extent, leaving some
space for reducing the regularization and hence allowing longer steps.

We note at this stage that the condition $s_k \in \calR_k$ implies that
\[
\calC_k \subseteq \calC_k^+
\tim{ and }
\calW_k^+ \subseteq \calW_k.
\]

Note that the above algorithm considerably simplifies in the Lipschitzian case
where $\calH = \emptyset$, since
\[
f_{\calW_k}(x) = f_\calM(x) = f(x)
\]
for all $k\geq 0$ and all $x \in \calF = \calF_\calQ$.

\numsection{Evaluation complexity for 'kernel-centered' fesible sets}
\label{complexity-s}

We start our worst-case analysis by formalizing our assumptions for problem
\req{psprob}.

\ass{AS.1}{ The feasible set $\calF$ is closed, convex and non-empty.}

\ass{AS.2}{ Each element function $f_i$ ($i \in \calN$) is $p$ times
  continuously differentiable in an open set containing $\calF$, where $p$ is
  odd whenever $\calH \neq \emptyset$.}

\ass{AS.3}{The $p$-th derivative of each $f_i$ ($i \in \calN$) is
  Lipschitz continuous on $\calF$ with associated Lipschitz constant $L_i$ (in
  the sense of \req{tensor-Lip-fi}).
}

\ass{AS.4}{There exists a constant $f_{\rm low}$ such that $f_\calN(x) \geq f_{\rm
    low}$ for all $x \in \calF$.
}

\ass{AS.5}{There exists a constant $\kappa^0_\calN \geq 0$ such that
\[
\|\nabla_x^jf_\calN(0)\| \leq \kappa_\calN^0
\]
for all $j \in \ii{p}$.
}

\noindent
Note that AS.4 is necessary for problem \req{psprob} to be well-defined. Also
observe that AS.5 guarantees the existence of a constant $\kappa_\calN \geq
\kappa_\calN^0$ such that
\beqn{dernice}
\|\nabla_x^1f_\calN(x))\| \leq \kappa_\calN
\tim{for all} x \in \{ x \in \calF \mid \|x\| \leq  1\}.
\eeqn
Obviously, AS.2 alone implies \req{dernice} (without the need
of assuming AS.5) if $\calF$ is bounded.

We first observe that our assumptions on the partially separable nature of the
objective function imply the following useful bounds.

\llem{psnorm-l}{
There exist constants $0 < \varsigma_{\min} \leq \varsigma_{\max}$ such that,
for all $s \in \Re^m$ and all $v \geq 1$ and for any subset $\calX \subseteq \calM$,
\beqn{s-sumsi}
\varsigma_{\min}^v \|s_\calX\|^v
\leq \sum_{i\in \calX} \|s_i\|^v
\leq |\calX| \,\varsigma_{\max}^v \,\|s_\calX\|^v,
\eeqn
where $s_\calX = P_{\spanset_{i\in\calX}\{U_i^T\}}(s)$.
}

\proof{
Assume that, for every $\varsigma > 0$ there exists a vector $s_\varsigma$
in $\spanset_{i\in\calX}\{U_i^T\}$ of norm 1 such
that $\max_{i \in \calX} \|U_i s_\varsigma\| < \varsigma \|s_\varsigma\|=
\varsigma$.  Then taking a sequence of $\{\varsigma_i\}$ converging to zero
and using the compactness of the unit sphere, we obtain that the sequence
$\{s_{\varsigma_i}\}$ has at least one limit point $s_0$ with $\|s_0\|=1$  such
that $\max_{i\in\calX} \|U_i s_0 \| = 0$, which is impossible since we
assumed that the intersection of the nullspaces of the $U_i$ is reduced to the
origin.  Thus our assumption is false and there is constant $\varsigma_{\min}
> 0$ such that, for every $s \in \spanset_{i\in\calX}\{U_i^T\}$,
\[
\max_{i \in \calX} \|s_i\|
= \max_{i \in \calX} \|U_is\|
\geq \varsigma_{\min}\|s\|.
\]
The first inequality of \req{s-sumsi} then follows from the fact that
\[
\sum_{i\in\calX} \|s_i\|^v
\geq \max_{i \in \calX} \|s_i\|^v
\geq \varsigma_{\min}^v \|s\|^v.
\]
We have also that
\[
\sum_{i\in\calX} \|s_i\|^v
\leq |\calX| \max_{i \in \calX} \|U_is\|^v
\leq |\calX| \max_{i \in \calX}\left(\|U_i\|\|s\|\right)^v,
\]
which yields the second inequality of \req{s-sumsi} with $\varsigma_{\max}=
\bigmax_{i \in \calX}\|U_i\|$.
}

\noindent
Taken for $v=1$ and $\calX = \calN$, this lemma states that
$\sum_{i\in\calN}\|\cdot\|$ is a norm on $\Re^n$ whose equivalence constants
with respect to the Euclidean one are $\varsigma_{\min}$ and
$|\calN|\,\varsigma_{\max}$. In most applications, these constants are very moderate
numbers.

\noindent
We now turn to the consequence of the Lipschitz continuity of $\nabla_x^p
f_i$ and define, for a given $k \geq 0$ and a given constant $\phi>0$
independent of $\epsilon$,
\beqn{Ok-def}
\calO_{k,\phi} \eqdef \{ i \in \calW_k^+\cap \calH
\mid \min[\, |x_{i,k}|, \,|x_{i,k}+s_{i,k}|\,] \geq \phi \}.
\eeqn
Note that
\[
\calO_{k,\phi} = \calH \setminus \Big[ \calC(x_k,\phi) \cup
\calC(x_k+s_k,\phi) \Big].
\]

\llem{Lip-th}{
Suppose that AS.2 and AS.3 hold. Then, for $k \geq 0$ and $L_{\max} \eqdef
\max_{i\in\calN}L_i$,
\beqn{fi-modi-N}
f_i(x_{i,k}+s_{i,k})
= m_i(x_{i,k},s_{i,k}) +
   \frac{1}{(p+1)!}\Big[ \tau_{i,k} (p+1)L_{\max} -  \sigma_{i,k} \Big]\|s_{i,k}\|^{p+1}
\tim{ with }| \tau_{i,k} | \leq 1,
\eeqn
for all $i \in \calN$.
If, in addition,  $\phi >0$ is given and independent of
$\epsilon$, then there exists a constant $L(\phi)$ independent of $\epsilon$
such that
\beqn{gf-gmod}
\| \nabla_x^1 f_{\calN\cup\calO_{k,\phi}}(x_k+s_k)
   -  \nabla_s^1 m_{\calN\cup\calO_{k,\phi}}(x_k,s_k) \|
\leq L(\phi) \|s_k\|^p.
\eeqn
}

\proof{
First note that, if $f_i$ has a Lipschitz continuous $p$-th derivative as a
function of $U_ix$, then \req{derjU} shows that it also has a Lipschitz continuous
$p$-th derivative as a function of $x$.  It is therefore enough to consider
the element functions as functions of $x_i=U_ix$.

AS.3 and \req{tensor-Lip-fi} imply that
\beqn{fi-Lip}
f_i(x_{i,k}+s_{i,k})
= T_{f_i,p}(x_{i,k},s_{i,k}) + \frac{\tau_{i,k}}{p!} L_{\max} \|s_{i,k}\|^{p+1}
\tim{ with }| \tau_{i,k} | \leq 1,
\eeqn
for each $i \in \calN$ (see \cite{BirgGardMartSantToin17} or
\cite[Section~2.2]{CartGoulToin16a}), and
\req{fi-modi-N} then follows from \req{model-fi}.

Consider now $i\in \calO_{k,\phi}$ and assume first that
$x_{i,k}> \phi$ and $x_{i,k}+s_{i,k}> \phi$. Then $f_i(x_i) = x_i^q$ is infinitely
differentiable on the interval $[x_{i,k}, x_{i,k}+s_{i,k}] \subset [\phi, \infty)$ and
the norm of its $(p+1)$-st derivative tensor is bounded above on this interval
by
\beqn{LH-def}
L_\calH(\phi)\eqdef\left|\prod_{\ell=0}^{p+1} (q-\ell)\right|\phi^{q-p-1}.
\eeqn
We then apply the same reasoning as above using the Taylor series expansion of
$x_i^q$ at $x_{i,k}$ and, because of the first line of \req{modiH}, deduce
that
\beqn{fi-modi-O}
f_i(x_{i,k}+s_{i,k})
= m_i(x_{i,k},s_{i,k}) +
   \frac{1}{(p+1)!} \tau_{i,k} (p+1)L_\calH(\phi) |s_{i,k}|^{p+1}
\tim{ with }| \tau_{i,k} | \leq 1,
\eeqn
and
\beqn{gfi-gmodi-O}
\| \nabla_x^1 f_i(x_{i,k}+s_{i,k}) -  \nabla_s^1 m_i(x_{i,k},s_{i,k}) \|
\leq L_\calH(\phi) |s_{i,k}|^p,
\eeqn
hold in this case (see \cite{BirgGardMartSantToin17}).
The argument is obviously similar if
$x_{i,k}<- \phi$ and $x_{i,k}+s_{i,k}<- \phi$, using symmetry and the
second line of \req{modiH}. Let us now consider the case where $x_{i,k} >
\phi$ and $x_{i,k}+s_{i,k} < -\phi$.  The expansion \req{combined-model}
then shows that we may reason as for $x_{i,k}<- \phi$ and $x_{i,k}+s_{i,k}<-
\phi$ using a Taylor expansion at $-x_i$ (which we know by symmetry)
and the third line of \req{modiH}.  The case where $x_{i,k} <
-\phi$ and $x_{i,k}+s_{i,k} > \phi$ is similar, using the fourth line of \req{modiH}.
As a consequence, \req{fi-modi-O} and \req {gfi-gmodi-O} hold for every $i\in
\calO_{k,\phi}$  with Lipschitz constant $L_\calH(\phi)$.
Moreover, using \req{s-sumsi} and the
definitions \req{LH-def},
\[
\sum_{i\in \calN\cup\calO_{k,\phi}} L_i\|s_i\|^{p+1}
\leq \max\left[L_{\max},L_\calH(\phi)\right] \sum_{i\in\calN\cup\calO_{k,\phi}}\|s_i\|^{p+1}
\]
from which 
\req{gf-gmod} may in turn be derived from \req{gfi-gmodi-O} and \req{s-sumsi} with
\beqn{Ls-def}
L(\phi)
\eqdef |\calM| \,\varsigma_{\max}^{p+1}\, \max\left[L_{\max},L_\calH(\phi)\right].
\eeqn

}

\noindent
Note that there is no dependence on $\phi$ in $L$ if $\calH = \emptyset$.

We now return to our statement that
\beqn{kerUiinF}
\ker(U_i)\cap\calF\neq\emptyset
\eeqn
may be assumed without loss
of generality for all $i \in \calH$. Indeed, assume that \req{kerUiinF} fails
for $j \in \calH$.  Then  $j \in \calO_{k, \xi_j}$ for all $k \geq
0$, where $\xi_j>0$ is the distance between $\ker(U_j)$ and $\calF$,
and we may transfer $j$ from $\calH$ to $\calN$ (possibly
modifying $L_{\max}$).

The definition of the model in \req{full-model} also implies a simple lower
bound on model decrease.

\llem{mdecr}{
For all $k\geq 0$,
\beqn{Dphi}
\delta T_k \geq \frac{1}{(p+1)!}\, \sigma_{\min} \sum_{i \in \calN}\|s_{i,k}\|^{p+1},
\eeqn
$s_k \neq 0$ and \req{rhok-def} is well-defined.
}

\proof{
  The bound directly follows from \req{deltaT-def}, the observation that the algorithm
  enforces $\delta m_k > 0$ and \req{sig-decr}.  Moreover, $\chi_m(x_k,0,\epsilon)
  = \chi_f(x_k,\epsilon) > \epsilon$. As a consequence, \req{mod-term} cannot hold for
  $s_k=0$ since termination would have then occured in Step~1 of
Algorithm~\ref{psarp-algo}. Hence at least one $\|s_{i,k}\|$ is strictly positive
because of \req{s-sumsi} and \req{Dphi} therefore implies that \req{rhok-def} is
well-defined.  }

\noindent
We now verify that the two-sided model \req{mod-hir} is an overestimate of the
function $|x|^q$ for all relevant $x_i$ and $s_i$.

\llem{a-overestimation-l}{
  Suppose that AS.2 holds. Then, for $i\in \calH$ and all $x_i,s_i \in \Re^n$
  with $x_i \neq 0 \neq x_i+s_i$, we have that
  \beqn{a-overest}
  |x_i+s_i|^q \leq  m_i(x_i,s_i).
  \eeqn
}

\proof{
  Since $i\in\calH$ by assumption, this implies that $\calH \neq \emptyset$,
  and thus, by AS.2, that $p$ is odd.
  From the mean-value theorem, we obtain that
  \beqn{a-Tayl-rest}
  \begin{array}{lcl}
  |x_i+s_i|^q
  & = & |x_i|^q+q \bigsum_{j=1}^p \bigfrac{1}{j!}\left[\bigprod_{\ell=1}^{j-1}(q-\ell)\right]
        |x_i|^{q-j} \mu(x_i,s_i)^j \\*[2ex]
  &   &  + \bigfrac{1}{(p+1)!}\left[\bigprod_{\ell=1}^p(q-\ell)\right]
        |U_iz|^{q-p-1} \mu(x_i, s_i)^{p+1}
  \end{array}
  \eeqn
  for some $z$ such that, using symmetry, $z \in [ x, x+s]$ if
  $(U_ix)(U_i(x+s))>0$ or $z \in [-x,x+s]$ otherwise. As a consequence,
  we have that
  \[
  |U_iz| \geq \min[ \, |x_i|, \, |x_i+s_i| \,] > 0.
  \]
  Remember now that $p$ is odd.  Then, using that $q \in (0,1)$, we have that
  \[
  \mu(x_i,s_i)^{p+1} \geq 0
  \tim{ and }
  \bigprod_{\ell=1}^p(q-\ell) < 0.
  \]
  The inequality
  \beqn{a-T-overest}
  |x_i+s_i|^q
  \leq  |x_i|^q + q \bigsum_{j=1}^p \bigfrac{1}{j!}\left[\bigprod_{\ell=1}^{j-1}(q-\ell)\right]
        |x_i|^{q-j} \mu(x_i,s_i)^j
  \eeqn
  therefore immediately follows from \req{a-Tayl-rest}, proving \req{a-overest}.
} 

\noindent
We next investigate the consequences of the model's definition \req{mod-hir}
when the singularity at the origin is approached and show that the two-sided
model has to remain large along the steps when $x_{i,k}$ is not too far from
the singularity.

\llem{a-grad-lemma}{
Suppose that $p \geq 1$ is odd, $q \in (0,1)$, $i\in\calH$, $|x_i|\in(\epsilon,1]$,
and $|x_i+s_i| \geq \epsilon$. Then
\beqn{a-gbig}
|\nabla_{s_i}^1 m_i(x_i,s_i)| > \half q \,|\nabla_{s_i}^1 m_i(x_i,0)|.
\eeqn
}

\proof{
  Following the argument in the proof of Lemma~\ref{Lip-th}, it is sufficient to
  consider that $x_i >0$ and $x_i+s_i>0$.
  From \req{modiH} (where $\mu(x_i,s_i) = s_i$), we have that
  \beqn{a-grad-mod-hi}
  \nabla_{s_i}^1 T_{x^q,p}(x_i,s_i)
  = q\sum_{j=1}^p \frac{1}{(j-1)!}\left[\prod_{\ell=1}^{j-1}(q-\ell)\right]x_i^{q-j}s_i^{j-1}.
  \eeqn
  Define $s_i = \beta x_i$. This gives that
  \req{a-grad-mod-hi} now reads
  \beqn{a-s-series}
  \nabla_{s_i}^1 T_{x^q,p}(x_i,\beta x_i)
  = q\sum_{j=1}^p \frac{1}{(j-1)!}\left[\prod_{\ell=1}^{j-1}(q-\ell)\right]x_i^{q-1}\beta^{j-1},
  \eeqn
  from which we deduce that
  \beqn{gm0}
  \nabla_{s_i}^1 m_i(x_i,0) = \nabla_{s_i}^1 T_{x^q,p}(x_i,0) = qx_i^{q-1}.
  \eeqn
  Suppose first that $s_i <0$, i.e. $\beta \in (-1,0)$,
  and observe that $s_i^{j-1} < 0$ exactly whenever $\bigprod_{\ell=1}^{j-1}(q-\ell)
  < 0$, and thus, using $x_i \leq 1$ and \req{gm0}, that
  \beqn{a-gradbig1}
  \nabla_{s_i}^1 m_i(x_i,s_i) > qx_i^{q-1} = \nabla_{s_i}^1 m_i(x_i,0)
  \tim{ for } \beta \in (-1,0).
  \eeqn
  Suppose now that $\beta \in (0,\third)$. Then \req{a-s-series} implies that
  \[
  \begin{array}{lcl}
  \nabla_{s_i}^1 T_{x^q,p}(x_i,\beta x_i)
  & \geq & qx_i^{q-1} -q\bigsum_{j=2}^p \left|\bigfrac{1}{(j-1)!}
          \left[\bigprod_{\ell=1}^{j-1}(q-\ell)\right]\right|x_i^{q-1}(\third)^{j-1}\\*[1ex]
  & = & qx_i^{q-1}\left(1
          -\bigsum_{j=2}^p \left|\frac{1}{(j-1)!}\left[\bigprod_{\ell=1}^{j-1}(q-\ell)\right]
            \right|(\third)^{j-1}\right).
  \end{array}
  \]
  Observe now that
  \beqn{Pik-bound}
  \left|\frac{1}{(j-1)!}\left[\bigprod_{\ell=1}^{j-1}(q-\ell)\right]\right|
  = \left|\prod_{\ell=1}^{j-1}\frac{q-\ell}{\ell}\right|
  \leq 1,
  \eeqn
  and therefore
  \[
  \begin{array}{lcl}
  \nabla_{s_i}^1 T_{x^q,p}(x_i,\beta x_i)
  & \geq & qx_i^{q-1}\left(1 - \bigsum_{j=2}^p (\third)^{j-1}\right) \\*[3ex]
  & > & qx_i^{q-1}\left(1 - \bigsum_{j=2}^\infty (\third)^{j-1}\right) \\*[3ex]
  & = & qx_i^{q-1}\left(1-\frac{\third}{1-\third}\right).
  \end{array}
  \]
  Using \req{gm0}, this implies that
  \beqn{a-gradbig-2}
  \nabla_{s_i}^1 T_{x^q,p}(x_i,\beta x_i) \geq \frac{1}{2} \nabla_{s_i}^1 T_{x^q,p}(x_i,0)
  \tim{ for } \beta \in [0,\third].
  \eeqn
  Suppose therefore that
  \beqn{a-betabiggish}
  \beta > \third.
  \eeqn
  We note that \req{a-s-series} gives that
  \[
  \nabla_{s_i}^1 T_{x^q,1}(x_i,s_i) = qx_i^{q-1}
  \tim{ and }
  \nabla_{s_i}^1 T_{x^q,t+2}(x_i,s_i)= \nabla_{s_i}^1 T_{x^q,t}(x_i,s_i) + qx_i^{q-1}h_t(\beta)
  \]
  for $t\in \ii{p-2}$ odd, where
  \beqn{ht-def}
  \begin{array}{lcl}
  h_t(\beta)
  & \eqdef & \bigfrac{1}{t!}\left[\bigprod_{\ell=1}^t(q-\ell)\right]\beta^t
  + \bigfrac{1}{(t+1)!}\left[\bigprod_{\ell=1}^{t+1}(q-\ell)\right]\beta^{t+1}\\*[3ex]
  & = & \bigfrac{1}{t!}\left[\bigprod_{\ell=1}^t(q-\ell)\right]\beta^t
                   \left(1+\bigfrac{q-(t+1)}{t+1}\beta\right).
  \end{array}
  \eeqn
  It is easy to verify that $h_t(\beta)$ has a root of multiplicity $t$ at zero
  and another root
  \[
  \beta_{0,t} = \frac{t+1}{t+1-q} \in \left(1\,, \, \frac{2}{2-q}\right),
  \]
  where the last inclusion follows from the fact that $q \in (0,1)$.  We also
  observe that $h_t(\beta)$ is a polynomial of even degree (since $t$ is
  odd). Thus
  \beqn{troot}
  h_t(\beta) \geq 0
  \tim{ for all } \beta \geq \frac{t+1}{t+1-q} \tim{ and } t\in \ii{p} {\rm ~odd}.
  \eeqn
  Now
  \beqn{ratsumht}
  \begin{array}{lcl}
  \bigfrac{\nabla_{s_i}^1 T_{x^q,p}(x_i,\beta x_i)}{qx_i^{q-1}}
  &  =  & \bigfrac{\nabla_{s_i}^1T_{x^q,p-2}(x_i,\beta x_i)}{qx_i^{q-1}}+h_{p-2}(\beta)\\*[1.5ex]
  &  =  & \bigfrac{\nabla_{s_i}^1 T_{x^q,1}(x_i,\beta x_i)}{qx_i^{q-1}}
           + \bigsum_{j=1, \,j {\rm ~odd}}^{p-2} h_j(\beta)\\*[1.5ex]
  &  =  &  1 + \bigsum_{\mystack{j=1, \,j {\rm ~odd}}{h_j(\beta)<0}}^{p-2} h_j(\beta) +
                \bigsum_{\mystack{j=1, \,j {\rm ~odd}}{h_j(\beta)\geq 0}}^{p-2} h_j(\beta)\\*[1.5ex]
  &  \geq &  1 + \bigsum_{\mystack{j=1, \,j {\rm ~odd}}{h_j(\beta)<0}}^{p-2} h_j(\beta)
  \end{array}
  \eeqn
  where we used \req{a-s-series} to derive the third equality.  Observe now that,
  because of \req{troot},
  \beqn{t0-def}
  \begin{array}{lcl}
  \{ j \in \ii{p-2} {\rm ~odd} \mid h_j(\beta)<0 \}
  & = & \left\{ j \in \ii{p-2} {\rm ~odd} \mid \beta < \frac{t+1}{t+1-q} \right\}\\*[1ex]
  & \eqdef & \{ j \in \{1,\ldots, t_0\} | \mid j {\rm ~odd}\}
  \end{array}
  \eeqn
  for some  odd integer $t_0 \in \ii{p-2}$.  Hence
  we deduce from \req{ht-def} and \req{ratsumht} that
  \beqn{alternating}
  \bigfrac{\nabla_{s_i}^1 T_{x^q,p}(x_i,\beta x_i)}{qx_i^{q-1}}
  \geq 1 + \sum_{j=1}^{t_0+1} \frac{1}{j!}\left[\prod_{\ell=1}^j(q-\ell)\right]\beta^j.
  \eeqn
  Moreover, since $h_t(\beta)<0$ for $t\in\ii{t_0}$ odd and observing that
  the second term in the first right-hand side of \req{ht-def} is always positive for $t$
  odd, we deduce that the terms in the summation of \req{alternating}
  alternate in sign.  We also note that they are decreasing in absolute value
  since
  \[
  \frac{1}{(t+1)!}\left|\prod_{\ell=1}^{t+1}(q-\ell)\right|\beta^{t +1}
  - \frac{1}{t!}\left|\prod_{\ell=1}^t(q-\ell)\right|\beta^t
  = \frac{1}{t!}\left|\prod_{\ell=1}^t(q-\ell)\right|\beta^t
   \left(\frac{t+1-q}{t+1}\,\beta -1\right)
  \]
  and \req{troot} ensures that the term in brackets in the right-hand side is
  always negative for $q \in (0,1)$ and $t \in \ii{t_0}$ odd.  Thus, keeping
  the first (most negative) term in \req{alternating}, we obtain that
  \beqn{a-gradbig-4}
  \nabla_{s_i}^1 T_{x^q,p}(x_i,\beta x_i)
  \geq  qx_i^{q-1}(1 +(q-1)\beta)
  \geq \frac{q}{2-q} \nabla_{s_i}^1 T_{x^q,p}(x_i,0)
   > \frac{q}{2}\nabla_{s_i}^1 T_{x^q,p}(x_i,0).
  \eeqn
  where we used \req{a-s-series} to
  deduce the second inequality. Combining \req{a-gradbig1}, \req{a-gradbig-2}
  and \req{a-gradbig-4} then yields that \req{a-gbig} holds for all $\beta \in
  (-1,\infty)$, which completes the proof since $s_i=\beta x_i$.
} 

\noindent
Our next step is to verify that the regularization parameters
$\{\sigma_{i,k}\}_{i \in \calN}$ cannot grow unbounded.

\llem{sigmax}{
Suppose that AS.2 and AS.3 hold.  Then, for all $i \in \calN$ and all $k \geq 0$,
\beqn{sigma-max}
\sigma_{i,k} \in [\sigma_{\min}, \sigma_{\max}],
\eeqn
where $\sigma_{\max} \eqdef \gamma_2 (p+1)L_{\max}$.
}

\proof{
Assume that, for some $i \in \calN$ and $k \geq 0$, $\sigma_{i,k} \geq (p+1)L_i$.
Then \req{fi-modi-N} gives that \req{up-cond} must fail, ensuring
\req{sigma-max} because of the mechanism of the algorithm.
}

\noindent
We next investigate the consequences of the model's definition \req{mod-hir}
when the singularity at the origin is approached.

\llem{a-big-grad-l}{
Suppose that AS.2 and AS.5 (and thus \req{dernice}) hold and that $\calH \neq
\emptyset$. Let
\beqn{a-omega-def}
\omega \eqdef \min\left[ 1,
\left(\frac{4\left[ p\,\kappa_\calN+\frac{|\calN|}{p!}\varsigma^p\sigma_{\max}\right]}
          {q^2}\right)^{\frac{1}{q-1}}
           \right],
\eeqn
and suppose, in addition, that
\beqn{a-sknleq1}
\|s_k\| \leq 1
\eeqn
and that, for some $i \in \calH$,
\beqn{a-assxik}
|x_{i,k}| \in (0,\omega).
\eeqn
Then
\beqn{a-big-grad}
\|P_{\calR{\{i\}}}[\nabla_s^1 m(x_k,s_k)]\| \geq \quarter q^2\omega^{q-1}
\tim{ and } \sgn\left(P_{\calR{\{i\}}}[\nabla_s^1 m(x_k,s_k)]\right) = \sgn(x_{i,k}+s_{i,k})
\eeqn
where $\calR_{\{i\}} \eqdef \spanset\{U_i^T\}$.
}

\proof{
Consider $i\in \calH$. Suppose, for the sake of simplicity, that
\beqn{a-simplicity}
x_{i,k} > 0 \tim{ and } x_{i,k}+s_{i,k} > 0.
\eeqn
We first observe that Lemma~\ref{a-grad-lemma}  implies that
\beqn{a-gmbig}
\nabla_{s_i}^1m_i(x_{i,k},s_i) \geq \half q \nabla_{s_i}^1m_i(x_{i,k},0)
\tim{ for all } s_i \neq -x_{i,k}.
\eeqn
Moreover,
\[
\nabla_s^1 m_\calN(x_k,s_k)
= \nabla_x^1 f_\calN(x_k)
  + \sum_{j=2}^p \frac{1}{(j-1)!}\nabla_x^j f_\calN(x_k)[s_k]^{j-1}
    + \frac{1}{p!}\sum_{\ell \in \calN}\sigma_{\ell,k}s_{\ell,k}\|s_{\ell,k}\|^{p-1}
\]
and thus, using the contractive property of
orthogonal projections,\req{a-sknleq1}, \req{dernice} and \req{s-sumsi}, that
\beqn{a-prevbound}
\begin{array}{lcl}
\|P_{\calR{\{i\}}}[\nabla_s^1 m_\calN(x_k,s_k)]\|
& \leq & \|\nabla_s^1 m_\calN(x_k,s_k)\| \\
& \leq & \kappa_\calN [1+(p-1)] + \bigfrac{|\calN|}{p!}\varsigma^p\sigma_{\max}\\
&   =  & p \, \kappa_\calN + \bigfrac{|\calN|}{p!}\varsigma^p\sigma_{\max}.
\end{array}
\eeqn
We next successively use the linearity of $P_{\calR{\{i\}}}[\cdot]$, the
triangle inequality, \req{a-gmbig}, the facts that
\[
\|U_i^T\|=1
\tim{ and }
|\nabla_{s_i}^1 m_i(x_k,s_k)|= q|x_{i,k}|^{q-1}\geq q\omega^{q-1},
\]
the bound \req{a-prevbound}, and \req{a-omega-def} to deduce that
\[
\begin{array}{lcl}
\|P_{\calR{\{i\}}}[\nabla_s^1 m(x_k,s_k)]\|
& = & \|P_{\calR{\{i\}}}[\nabla_s^1 m_\calN(x_k,s_k)+\nabla_s^1\bigsum_{j\in\calH} m_j(x_k,s_k)]\|\\*[1ex]
& = & \|P_{\calR{\{i\}}}[\nabla_s^1 m_\calN(x_k,s_k)+\bigsum_{j\in\calH}U_j^T\nabla_{s_j}^1 m_j(x_k,s_k)]\|\\*[1ex]
& = & \|P_{\calR{\{i\}}}[\nabla_s^1 m_\calN(x_k,s_k)]+U_i^T\nabla_{s_i}^1 m_i(x_k,s_k)\|\\*[1ex]
& \geq & \Big| \|U_i^T\nabla_{s_i}^1 m_i(x_k,s_k)\|
               -\|P_{\calR{\{i\}}}[\nabla_s^1 m_\calN(x_k,s_k)]\| \Big| \\*[1ex]
& \geq & \half q^2\omega^{q-1}-\left[p\,\kappa_\calN
                 +\bigfrac{|\calN|}{p!}\varsigma^p\sigma_{\max}\right]\\*[1ex]
& \geq & \quarter q^2 \omega^{q-1},
\end{array}
\]
which proves the first part of \req{a-big-grad} and, because of \req{a-gmbig},
implies the second, for the case where \req{a-simplicity} holds.  The
proof for the cases where
\[
\Big[\, x_{i,k} < 0 \tim{ and } x_{i,k}+s_{i,k} < 0 \, \Big]
\tim{ or } x_{i,k}(x_{i,k}+s_{i,k}) < 0
\]
are identical when making use of the symmetry $m_i(x_i)$ with respect to the
origin.
}  

\noindent
Note that, like $\sigma_{\max}$, $\omega$ and $\beta$ only depend on
problem data. In particular, they are independent of $\epsilon$.
Lemma~\ref{a-big-grad-l} has the following crucial consequence.

\llem{a-away-l}{
Suppose that AS.2, AS.5 and the assumptions \req{a-sknleq1}--\req{a-assxik} of
Lemma~\ref{a-big-grad-l} hold and that $\calH \neq \emptyset$. Suppose in
addition that \req{mod-term} holds at $x_k,s_k$. Then, either
\beqn{a-xplus-away}
|x_{i,k}+s_{i,k} | \leq \epsilon
\tim{ or }
|x_{i,k}+s_{i,k} | \geq \omega
\ms
(i \in \calH).
\eeqn
}

\proof{
If $j \in \calH \cap \calC(x_k+s_k, \epsilon)$, then clearly
$|x_{j,k}+s_{j,k} | \leq \epsilon$, and there is nothing more to prove.
Consider therefore any
$j \in\calH \setminus \calC_k^+\subseteq \calW_k^+$
and observe that the separable nature of the linear optimization problem
in \req{chim-def} implies that
\beqn{a-pe1}
\begin{array}{lcl}
\left|\bigmin_{\mystack{x_k+s_k+d \in \calF}{d \in \calR_{\{j\}},\|d\|\leq 1}}
      P_{\calR_{\{j\}}}[\nabla_s^1m(x_k,s_k)]^Td
      \right|
& = & \left|\bigmin_{\mystack{x_k+s_k+d \in \calF}{d \in \calR_{\{j\}},\|d\|\leq 1}}
      \nabla_s^1m_{\calW_k^+}(x_k,s_k)^Td
      \right| \\*[4.5ex]
& \leq & \left|\bigmin_{\mystack{x_k+s_k+d \in \calF}{ d \in \calR_k^+,\|d\|\leq 1}}
      \nabla_s^1 m_{\calW_k^+}(x_k,s_k)^Td
      \right| \\*[4.5ex]
& = & \chi_m(x_k,s_k,\epsilon)\\*[1ex]
& \leq & \quarter q^2 |x_{j,k}+s_{j,k}|^r.
\end{array}
\eeqn
Observe now that, because of the second part of \req{a-big-grad}
and the fact that $n_j=1$ because of \req{UiH-conds}, the optimal value for
the convex optimization problem in the left-hand side of this relation is
given by
\[
|P_{\calR_{\{j\}}}[\nabla_s^1m(x_k,s_k)]| \, |d_*|
\]
where $d_*$ is the problem solution and $d_*$ has the opposite sign of
$P_{\calR_{\{j\}}}[\nabla_s^1m(x_k,s_k)]$. Moreover, the facts that $j \in
\calH$ and \req{Fbox}  ensure that
$x_{j,k}+s_{j,k}+d_j = 0$ is feasible for the optimization problem on the
left-hand side of \req{a-pe1}, and hence that $|d_*| \geq |x_{j,k}+s_{j,k}|$.
Hence, we obtain that
\[
\quarter q^2\omega^{q-1} |x_{j,k}+s_{j,k}| \leq \quarter q^2  |x_{j,k}+s_{j,k}|^r,
\]
and thus, since $\omega\leq 1$, that
\[
|x_{j,k}+s_{j,k}| \geq \omega^{\frac{q-1}{r-1}} \geq \omega,
\]
and the second alternative in \req{a-xplus-away} holds.
} 

\noindent
The rest of our complexity analysis depends on the following partitioning of
the set of iterations. Let the index set of the ``successful'' and
``unsuccessful'' iterations be given by
\[
\calS \eqdef \{ k \geq 0 \mid \rho_k \geq \eta \}
\tim{and}
\calU \eqdef \{ k \geq 0 \mid \rho_k < \eta \}.
\]
We next focus on the case where $\calH \neq \emptyset$ and partition $\calS$
into subsets depending on $|x_{i,k}|$ and $|x_{i,k}+s_{i,k}|$ for $i \in
\calH$.  We first isolate the set of sucessful iterations which ``deactivate'' some
variable, that is
\[
\calS_{\epsilon}
\eqdef \{ k \in \calS \mid |x_{i,k}+s_{i,k}| \leq \epsilon \tim{for some} i \in \calH\},
\]
as well as the set of successful iterations with large steps
\beqn{a-Sn-def}
\calS_{\|s\|}
\eqdef \{ k \in \calS \setminus \calS_\epsilon \mid \|s_k\| > 1 \}.
\eeqn
Let us now choose a constant $\alpha \geq 0$ such that
\beqn{a-alpha-def}
\alpha  = \left\{ \begin{array}{ll}
 \threequarters \omega & \tim{if } \calH  \neq \emptyset, \\
  0                    & \tim{otherwise.}
\end{array}\right.
\eeqn
Then, at iteration $k \in \calS \setminus ( \calS_\epsilon \cup
\calS_{\|s\|})$, we distinguish
\[
\begin{array}{ll}
\calI_{\heartsuit,k} \eqdef  \Big\{ i \in \calH\setminus\calC_k \mid &
|x_{i,k}| \in [\alpha, +\infty)
\tim{and} |x_{i,k}+s_{i,k}| \in [\alpha, +\infty) \Big\}, \\*[1ex]
\calI_{\diamondsuit,k} \eqdef \Big\{ i \in \calH \setminus\calC_k \mid &
\Big( |x_{i,k}| \in [\omega,+\infty)
\tim{and} |x_{i,k}+s_{i,k}| \in (\epsilon,\alpha) \Big)\\*[1ex]
&  \tim{or}
\Big( |x_{i,k}| \in (\epsilon,\alpha)
\tim{and} |x_{i,k}+s_{i,k}| \in [\omega,+\infty)  \Big) \Big\},\\*[1ex]
\calI_{\clubsuit,k} \eqdef \Big\{ i \in \calH \setminus\calC_k \mid &
|x_{i,k}| \in (\epsilon, \omega)
\tim{and} |x_{i,k}+s_{i,k}| \in (\epsilon, \omega)\Big\}. \\*[1ex]
\end{array}
\]
Using these notations, we further define
\[
\begin{array}{ll}
  \calS_\heartsuit   \eqdef \{ k \in \calS \setminus (\calS_\epsilon \cup \calS_{\|s\|}) \mid
                             \calI_{\heartsuit,k} = \calH\setminus\calC_k \}, &
  \calS_\diamondsuit  \eqdef \{ k \in \calS \setminus (\calS_\epsilon \cup \calS_{\|s\|}) \mid
                             \calI_{\diamondsuit,k} \neq \emptyset \}, \\
  \calS_\clubsuit    \eqdef \{ k \in \calS \setminus (\calS_\epsilon \cup \calS_{\|s\|}) \mid
                             \calI_{\clubsuit,k} \neq \emptyset \}.
\end{array}
\]
Figure~\req{a-kind-of-steps-fig} displays the various kinds of steps
in $\calS_{\heartsuit,k}$, $\calS_{\diamondsuit,k}$, $\calS_{\clubsuit,k}$
and $\calS_{\epsilon,k}$.
\begin{figure}[htbp]
\begin{center}
\vspace*{1.5mm}
\includegraphics[height=6.5cm]{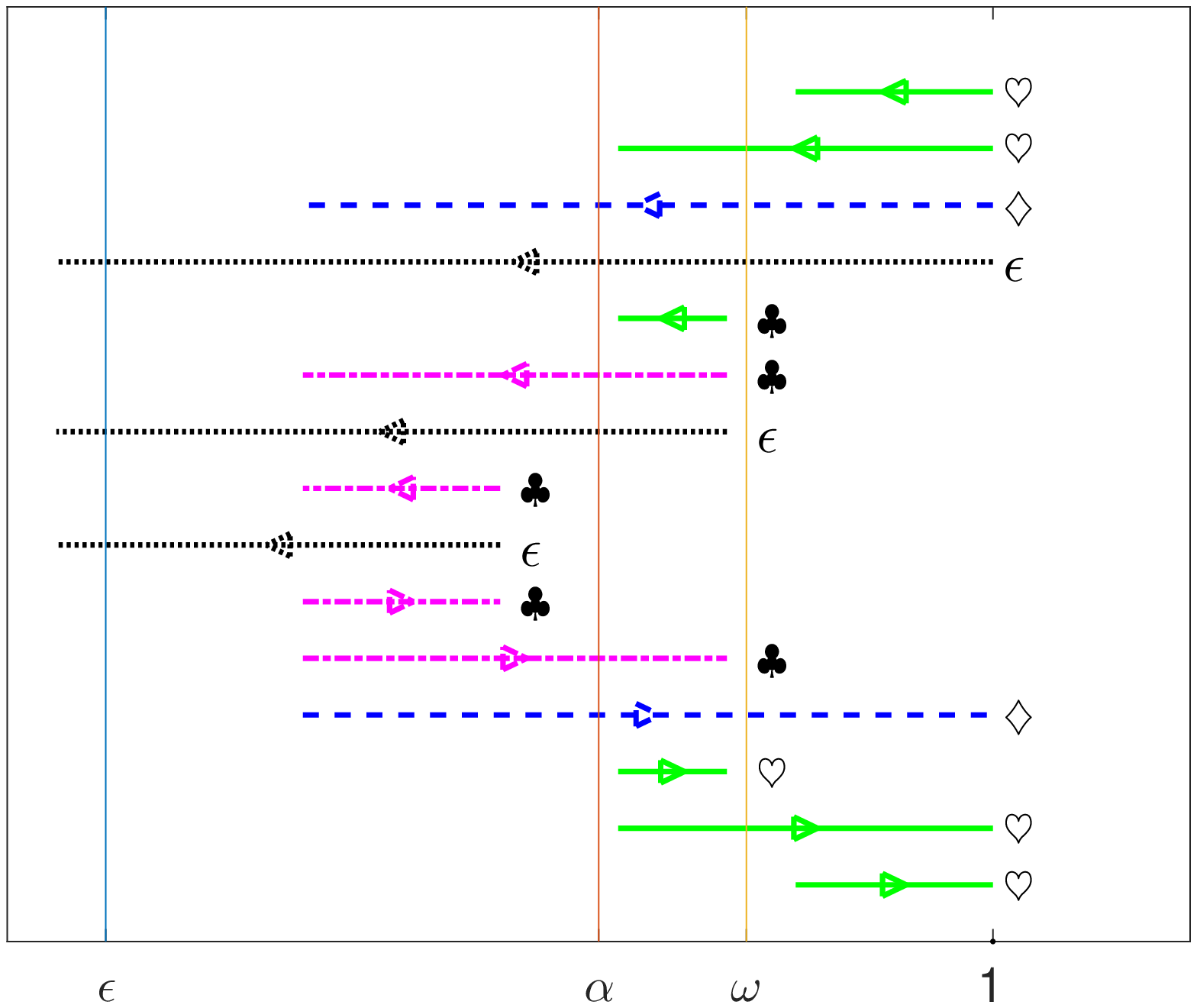} 
\caption{\label{a-kind-of-steps-fig} The various steps in $\calS\setminus\calS_{\|s\|}$
  depending on intervals containing their origin $|x_{i,k}|$ and end $|x_{i,k}+s_{i,k}|$ points.
  The vertical lines show, in increasing order, $\epsilon$, $\alpha$ and $\omega$. The
  line type of the represented step indicates that it belongs to
  $\calS_{\epsilon,k}$ (dotted), $\calS_{\heartsuit,k}$ (solid),
  $\calS_{\diamondsuit,k}$ (dashed) and $\calS_{\clubsuit,k}$
  (dash-dotted).
  The vertical axis is meaningless.}
\end{center}
\vspace*{-5mm}
\end{figure}

It is important to observe that the mechanism of the algorithm ensures that,
once an $x_i$ falls in the interval $[-\epsilon,\epsilon]$ at iteration $k$, it never
leaves it (and essentially ``drops out'' of the calculation).  Thus there are no
right-oriented dotted steps in Figure~\ref{a-kind-of-steps-fig} and also
\beqn{a-S-eps}
|\calS_\epsilon| \leq |\calH|.
\eeqn
Moreover Lemma~\ref{a-away-l} ensures that $\calI_{\clubsuit,k} = \emptyset$
for all $k \in \calS$, and hence that
\beqn{a-S-club}
|\calS_\clubsuit| = 0.
\eeqn
As a consequence, one has that $\calS_\epsilon$, $\calS_{\|s\|}$, $\calS_\heartsuit$,
and $\calS_\diamondsuit$ 
form a partition of $\calS$.
It is also easy to verify that, if $k \in \calS_{\diamondsuit}$ and $i \in
\calI_{\diamondsuit,k}$, then
\beqn{a-s-S-diam}
\|s_k\| \geq \|P_{\calR_{\{i\}}}(s_k)\| = |s_{i,k}| \geq \omega-\alpha = \quarter
\omega > 0,
\eeqn
where we have used the contractive property of orthogonal projections.

We now show that the steps at iterations whose index is in $\calS_\heartsuit$
are not too short.

\llem{a-schi+-lemma}{
Suppose that AS.1-AS.3 and AS.5 hold, that
\beqn{eps-small-1}
\epsilon <  \alpha
\eeqn
and consider $k \in \calS_\heartsuit$ before termination.  Then
\beqn{a-schi+}
\|s_k\| \geq (\kappa_\heartsuit \, \epsilon)^{\frac{1}{p}},
\eeqn
where
\beqn{a-kappaheart-def}
\kappa_\heartsuit
\eqdef \left[2(L(\alpha) +\theta+\bigfrac{|\calN|}{p!}\,\varsigma_{\max}^{p+1}\,\sigma_{\max})\right]^{-1}.
\eeqn
}

\proof{
Observe first that, since $k \in \calS_\heartsuit\subseteq \calS$, we have that $x_{k+1} =
x_k+s_k$ and, because $\epsilon \leq \alpha$ and $\calC_k^+\subseteq \calC_k$,
we deduce that $\calC_k=\calC_k^+=\calC_{k+1}$ and $\calR_k= \calR_k^+ =
\calR_{k+1}$.  Moreover the definition of $\calS_\heartsuit$ ensures that, for
all $i \in \calH\setminus\calC_k$,
\beqn{heartguy}
\min\Big[|x_{i,k}|,|x_{i,k}+s_{i,k}|\Big] \geq \alpha.
\eeqn
Hence
\[
\calO_*
\eqdef \calO_{k, \alpha}
= \calH\setminus\calC_k=\calH\setminus\calC_k^+,
\]
and thus
\beqn{a-Wstar}
\calR_* \eqdef \calR_k = \calR_k^+
\tim{ and }
\calW_* \eqdef \calW_k = \calW_k^+ = \calN \cup \calO_*.
\eeqn
As a consequence the step computation must have been completed because
\req{mod-term} holds, which implies that
\beqn{a-chim-ok}
\chi_m(x_k,s_k,\epsilon) =\chi_{m_{\calW_*}}(x_k,s_k,\epsilon) =
\left| \min_{\mystack{x_k+s_k+d \in \calF}{d\in\calR_*, \|d\|\leq 1}}
      \nabla_sm_{\calW_*}(x_k,s_k)^Td \right|
\leq \theta \|s_k\|^p.
\eeqn
Observe also that \req{a-Wstar}, \req{gf-gmod} with $\phi=\alpha$ (because $k
\in \calS_\heartsuit$) , \req{sigma-max} and \req{s-sumsi} then imply that
\beqn{a-errdqm}
\begin{array}{lcl}
\|\nabla_x^1 f_{\calW_*}(x_{k+1}) - \nabla_s^1 m_{\calW_*}(x_k,s_k)\|
&   =  & \|\nabla_x^1 f_{\calN\cup\calO_*}(x_{k+1})
         - \nabla_s^1 m_{\calN\cup\calO_*}(x_k,s_k)\| \\*[3ex]
& \leq & L(\alpha) \|s_k\|^p +\bigfrac{1}{(p+1)!}\, \sigma_{\max}
         \bigsum_{i \in \calN} \| \, \nabla_s^1 \|s_{i,k}\|^{p+1} \,\| \\*[3ex]
& \leq & L(\alpha) \|s_k\|^p
         +\bigfrac{1}{p!}\, \sigma_{\max}\bigsum_{i \in \calN} \|s_{i,k}\|^p \\*[3ex]
& \leq & L(\alpha) \|s_k\|^p
         + \bigfrac{|\calN|}{p!}\,\varsigma_{\max}^{p+1}\, \sigma_{\max}\|s_k\|^p \\*[3ex]
& = & \left[L(\alpha)+\bigfrac{|\calN|}{p!}\,\varsigma_{\max}^{p+1}\,\sigma_{\max}\right]
         \|s_k\|^p,
\end{array}
\eeqn
and also that
\beqn{a-con2-ARC2CC-sl-1}
\begin{array}{lcl}
\chi_f(x_{k+1},\epsilon)
&   =    & |\nabla_x^1 f_{\calW_*}(x_{k+1})[d_{k+1}]| \\*[1ex]
&  \leq  & |\nabla_x^1 f_{\calW_*}(x_{k+1})[d_{k+1}]
             - \nabla_s^1 m_{\calW_*}(x_k,s_k)[d_{k+1}]| \\*[1ex]
&        & \hspace*{1cm}  + |\nabla_s^1 m_{\calW_*}(x_k,s_k)[d_{k+1}]|,
\end{array}
\eeqn
where the first equality defines the vector $d_{k+1}$  with
\beqn{a-dplusbound-chi}
\|d_{k+1}\| \leq 1.
\eeqn
Assume now, for the purpose of deriving a contradiction, that
\beqn{a-schi+chi}
\|s_k\| <  \left[ \frac{\chi_f(x_{k+1},\epsilon)}
  {2(L(\alpha)+\theta+\bigfrac{|\calN|}{p!}\,\varsigma_{\max}^{p+1}\,\sigma_{\max})}\right]^{\frac{1}{p}}
\eeqn
at iteration $k \in \calS_\heartsuit$. Using \req{a-dplusbound-chi} and
\req{a-errdqm}, we then obtain that
\beqn{a-dgradj}
\begin{array}{lcl}
  \lefteqn{-\nabla_x^1 f_{\calW_*}(x_{k+1})[d_{k+1}]
            +\nabla_s^1 m_{\calW_*}(x_k,s_k)[d_{k+1}]}\\*[1.5ex]
\hspace*{4cm} &
\leq &|\nabla_x^1 f_{\calW_*}(x_{k+1})[d_{k+1}]
        - \nabla_s^1 m_{\calW_*}(x_k,s_k)[d_{k+1}] | \\*[1.5ex]
&  =  &|(\nabla_x^1 f_{\calW_*}(x_{k+1})
        - \nabla_s^1 m_{\calW_*}(x_k,s_k))[d_{k+1}] | \\*[1.5ex]
&\leq &\|\nabla_x^1 f_{\calW_*}(x_{k+1})
        - \nabla_s^1 m_{\calW_*}(x_k,s_k)\| \; \|d_{k+1}\| \\*[1.5ex]
&  <  & ( L(\alpha)+\bigfrac{|\calN|}{p!}\,\varsigma_{\max}^{p+1}\,\sigma_{\max}) \|s_k\|^p.
\end{array}
\eeqn
From \req{a-schi+chi} and the first part of \req{a-con2-ARC2CC-sl-1}, we have that
\[
\begin{array}{lcl}
-\nabla_x^1 f_{\calW_*}(x_{k+1})[d_{k+1}]
   + \nabla_s^1 m_{\calW_*}(x_k,s_k))[d_{k+1}]
& < &\half \chi_f(x_{k+1},\epsilon)\\*[1ex]
& = &  - \half \nabla_x^1 f_{\calW_*}(x_{k+1})[d_{k+1}],
\end{array}
\]
which in turn ensures that
\[
\nabla_s^1 m_{\calW_*}(x_k,s_k)[d_{k+1}]
< \half \nabla_x^1 f_{\calW_*}(x_{k+1})[d_{k+1}] < 0.
\]
Moreover, by definition of $\chi_f(x_{k+1},\epsilon)$,
\[
x_{k+1}+d_{k+1} \in \calF
\tim{and}
d_{k+1} \in \calR_{k+1} = \calR_k^+.
\]
Hence, using \req{chim-def} and \req{a-dplusbound-chi},
\beqn{a-nmdleqchi-chi}
|\nabla_s^1 m_{\calW_*}(x_k,s_k)[d_{k+1}]| \leq \chi_{m_{\calW_*}}(x_k,s_k,\epsilon).
\eeqn
We may then substitute this inequality in \req{a-con2-ARC2CC-sl-1} to deduce as
above that
\beqn{a-con2-ARC2CC-sl-1b}
\begin{array}{lcl}
\chi_f(x_{k+1})
& \leq & |\nabla_x^1 f_{\calW_*}(x_{k+1})[d_{k+1}]
          - \nabla_s^1 m_{\calW_*}(x_k,s_k)[d_{k+1}]|
         + \chi_{m_{\calW_*}}(x_k,s_k,\epsilon) \\*[2ex]
& \leq & (L(\alpha) + \theta + \bigfrac{|\calN|}{p!}\,\varsigma_{\max}^{p+1}\,\sigma_{\max})
         \|s_k\|^p
\end{array}
\eeqn
where the last inequality results from \req{a-dgradj}, the identity $x_{k+1} =
x_k+s_k$ and \req{a-chim-ok}.  But this contradicts our assumption that
\req{a-schi+chi} holds. Hence \req{a-schi+chi} must fail. The inequality \req{a-schi+}
then follows by combining this conclusion with the fact that
$\chi_f(x_{k+1},\epsilon) > \epsilon$ before termination.
}

\noindent
We are now ready to consider our first complexity result, whose proof uses
restrictions of the successful and unsuccessful iteration index sets defined
above to $\iibe{0}{k}$, which are given by
\beqn{a-Sk-Uk-def}
\calS_k \eqdef \iibe{0}{k} \cap \calS,
\ms
\calU_k \eqdef \iibe{0}{k}\setminus \calS_k,
\eeqn
respectively.

\lthm{a-comp1}{
Suppose that AS.1-AS.5 hold and that
\beqn{a-eps-upper}
\epsilon \leq  \left[ \alpha,
  \left(\quarter \omega\kappa_\heartsuit^{-\frac{1}{p+1}}\right)^p\right]
\tim{ if } \calH \neq \emptyset.
\eeqn
Then Algorithm~\ref{psarp-algo} requires at most
\beqn{a-succ-compl}
\kappa_\calS(f(x_0)-f_{\rm low})
   \epsilon^{-\sfrac{p+1}{p}}
  +|\calH|
\eeqn
successful iterations to return a point $x_\epsilon \in \calF$ such that
$\chi_f(x_\epsilon,\epsilon) \leq \epsilon$, for
\beqn{a-kapS-def}
\kappa_\calS = \bigfrac
      {(p+1)!}{\eta \, \sigma_{\min}\, \varsigma_{\min}^{p+1} }
      \Big[2( L(\alpha)+\theta
        +\bigfrac{|\calN|}{p!}\,\varsigma_{\max}^{p+1}\,\gamma_2) \Big]^{\frac{p+1}{p}}.
 \eeqn
}

\proof{
Let $k \in \calS$ be index of a successful iteration before termination, and
suppose first that $\calH \neq \emptyset$.
Because the iteration is successful, we obtain, using AS.4 and
Lemma~\ref{mdecr}, that
\beqn{a-lowb1}
f(x_0)-f_{\rm low}
\geq f(x_0)-f(x_{k+1})
\geq \bigsum_{\ell \in \calS_k}\Big[f(x_\ell)-f(x_\ell+s_\ell)\Big]
\geq \eta \bigsum_{\ell \in \calS_k} \Big[f(x_\ell)-T_{f,p}(x_\ell,s_\ell)\Big].
\eeqn
In addition to \req{a-Sk-Uk-def}, let us define
\beqn{a-Snk-Sepsk-def}
\calS_{\epsilon,k} \eqdef \iibe{0}{k} \cap \calS_\epsilon,
\ms
\calS_{\|s\|,k} \eqdef \iibe{0}{k} \cap \calS_{\|s\|},
\eeqn
\[
\calS_{\heartsuit,k} \eqdef \iibe{0}{k} \cap \calS_\heartsuit,
\ms
\calS_{\diamondsuit,k} \eqdef \iibe{0}{k}\cap \calS_\diamondsuit.
\]
We now use the fact that
$\calS_{\|s\|,k}\cup\calS_{\heartsuit,k}\cup\calS_{\diamondsuit,k}
= \calS_k\setminus\calS_{\epsilon,k}\subseteq\calS_k$,
and \req{s-sumsi} to deduce from \req{a-lowb1} that
\[
\begin{array}{lcl}
f(x_0)-f_{\rm low}
& \geq & \eta \left\{
              \bigsum_{\ell \in \calS_{\|s\|,k}}
                  \Big[f(x_\ell)-T_{f,p}(x_\ell,s_\ell)\Big] +
              \bigsum_{\ell \in \calS_{\heartsuit,k}}
              \Big[f(x_\ell)-T_{f,p}(x_\ell,s_\ell)\Big] \right.\\*[4ex]
&       & \hspace*{1cm} + \left.
              \bigsum_{\ell \in \calS_{\diamondsuit,k}}
                  \Big[f(x_\ell)-T_{f,p}(x_\ell,s_\ell)\Big]
              \right\}\\*[4 ex]
& \geq & \bigfrac{\eta\sigma_{\min}}{(p+1)!}
              \left\{
              |\calS_{\|s\|,k}|\bigmin_{\ell \in \calS_{\|s\|,k}}
                   \left[\bigsum_{i\in\calN}\|s_{i,\ell}\|^{p+1}\right]+
              |\calS_{\heartsuit,k}|\bigmin_{\ell \in \calS_{\heartsuit,k}}
                   \left[\bigsum_{i\in\calN}\|s_{i,\ell}\|^{p+1}\right]\right.\\*[3ex]
&       & \hspace*{1cm} + \left.
              |\calS_{\diamondsuit,k}|\bigmin_{\ell \in \calS_{\diamondsuit,k}}
                   \left[\bigsum_{i\in\calN}\|s_{i,\ell}\|^{p+1}\right]
            \right\}\\*[3 ex]
& \geq & \bigfrac{\eta\sigma_{\min}\varsigma_{\min}^{p+1}}{(p+1)!}
            \left\{
      |\calS_{\|s\|,k}|\bigmin_{\ell \in \calS_{\|s\|,k}} \|s_\ell\|^{p+1}+
      |\calS_{\heartsuit,k}|\bigmin_{\ell \in \calS_{\heartsuit,k}} \|s_\ell\|^{p+1}
             \right.\\*[3ex]
&       & \hspace*{1cm} + \left.
      |\calS_{\diamondsuit,k}|\bigmin_{\ell \in \calS_{\diamondsuit,k}} \|s_\ell\|^{p+1}
            \right\}.
\end{array}
\]
Because of of \req{a-Sn-def}, \req{a-Snk-Sepsk-def}, Lemma~\ref{a-schi+-lemma} and
\req{a-s-S-diam}, this now yields that
\[
\begin{array}{lcl}
f(x_0)-f_{\rm low}
& \geq & \bigfrac{\eta\sigma_{\min}\varsigma_{\min}^{p+1}}{(p+1)!}
            \left\{
            |\calS_{\|s\|,k}|+
            |\calS_{\heartsuit,k}| (\kappa_\heartsuit\epsilon)^{\frac{p+1}{p}}+
            |\calS_{\diamondsuit|,k} (\omega-\alpha)^{p+1}
            \right\} \\*[3 ex]
& \geq & \bigfrac{\eta\sigma_{\min}\varsigma_{\min}^{p+1}}{(p+1)!}
            \Big\{ |\calS_{\|s\|,k}|+|\calS_{\heartsuit,k}|
                 + |\calS_{\diamondsuit,k}|\Big\}
        \min\left[ (\kappa_\heartsuit\epsilon)^{\frac{p+1}{p}},(\quarter
          \omega)^{p+1}\right]\\*[3 ex]
& \geq & \bigfrac{\eta\sigma_{\min}\varsigma_{\min}^{p+1}}{(p+1)!}
        \,|\calS_k\setminus\calS_\epsilon|\,(\kappa_\heartsuit\epsilon)^{\frac{p+1}{p}}
\end{array}
\]
where we used \req{a-eps-upper}, the partition of $\calS_k\setminus
\calS_{\epsilon,k}$ in $\calS_{\|s\|,k}\cup\calS_{\heartsuit,k}
\cup\calS_{\diamondsuit,k}$ and the inequality
$\quarter \omega< 1$ to obtain the last inequality. Thus
\beqn{a-part-comp}
|\calS_k| \leq \kappa_\calS (f(x_0)-f_{\rm low}) \epsilon^{-\frac{p+1}{p}} + |\calS_{\epsilon,k}|,
\eeqn
where $\kappa_\calS$ is given by \req{a-kapS-def}.  The desired iteration
complexity \req{a-succ-compl} then follows from this bound,
$|\calS_{\epsilon,k}| \leq |\calS_\epsilon|$ and \req{a-S-eps}.
}

\noindent
To complete our analysis in terms of evaluations rather than successful
iterations, we need to bound the total number of all (successful and
unsuccessful) iterations.

\llem{a-succ-unsucc}{
Assume that AS.2 and AS.3 hold.  Then, for all $k \geq 0$,
\[
k \leq \kappa^a |\calS_k| + \kappa^b,
\]
where
\[
\kappa^a \eqdef  1 + \frac{|\calN|\,|\log \gamma_0|}{\log \gamma_1}
\tim{ and }
\kappa^b \eqdef \frac{|\calN|}{\log \gamma_1}
               \log\left(\frac{\sigma_{\max}}{\sigma_{\min}}\right).
\]
}

\proof{
 For $i\in\calN$, define
 \[
 \calJ_{i,k} \eqdef \{ j \in \iibe{0}{k} \mid \mbox{\req{sig-incr} holds with }
 k \leftarrow j\},
 \]
 (the set of iterations where $\sigma_{i,j}$ is increased) and
 \[
 \calD_{i,k} \eqdef \{ j \in \iibe{0}{k} \mid \mbox{\req{sig-decr} holds with }
 k \leftarrow j\}
 \subseteq \calS_k
 \]
 (the set of iterations where $\sigma_{i,j}$ in decreased), the final inclusion
 resulting from the condition that $\rho_k\geq \eta$ in both
 \req{muchover-neg} and \req{muchover-pos}.
 Observe also that the mechanism of the algorithm, the fact that $\gamma_0 \in
 (0,1)$ and Lemma~\ref{sigmax} impose that, for each $i \in \calN$,
 \[
 \sigma_{\min}\gamma_1^{|\calJ_{i,k}|}\gamma_0^{|\calS_k|}
 \leq \sigma_{i,0}\gamma_1^{|\calJ_{i,k}|}\gamma_0^{|\calD_{i,k}|}
 \leq \sigma_{i,k}
 \leq \sigma_{\max}.
 \]
 Dividing by $\sigma_{\min}>0$ and taking logarithms yields that, for all $i \in
 \calN$ and all $k>0$,
 \beqn{cardUik}
 |\calJ_{i,k}|\log \gamma_1 + |\calS_k|\log \gamma_0
 \leq \log\left(\frac{\sigma_{\max}}{\sigma_{\min}}\right).
 \eeqn
 Note now that, if \req{up-cond} fails for all $i \in \calN$ and given that
 Lemma~\ref{a-overestimation-l} ensures that $f_i(x_i+s_i) \leq m_i(x_i,s_i)$ for
 $i\in\calH\setminus\calC_k^+$, then
 \[
 \delta f_k
 = \sum_{i\in\calW_k^+} \delta f_{i,k}
 \geq  \sum_{i\in\calW_k^+} \delta m_{i,k}
 = \delta m_k.
 \]
 Thus, in view of \req{rhok-def}, we have that $\rho_k \geq 1 > \eta$ and  iteration $k$
 is successful.  Thus, if iteration $k$ is unsuccessful, $\sigma_{i,k}$ is
 increased with \req{sig-incr} for at least one $i\in\calN$. Hence we deduce
 that
 \beqn{cardUk}
 |\calU_k|
 \leq \sum_{i\in\calN} |\calJ_{i,k}|
 \leq |\calN|\, \max_{i\in\calN}|\calJ_{i,k}|.
 \eeqn
 The desired bound follows from \req{cardUik} and \req{cardUk} by
 using the  fact that $k = |\calS_k| + |\calU_k| - 1 \leq |\calS_k| +
 |\calU_k|$, the  term -1 in the equality accounting for iteration 0.
}

\noindent
We may now state our main evaluation complexity result.

\lthm{final-comp}{
Suppose that AS.1, \req{Fbox}, AS.2-AS.5 and \req{a-eps-upper} hold. Then
Algorithm~\ref{psarp-algo} using models \req{mod-hir} for $i\in \calH$ requires at most
\beqn{a-better-comp}
\kappa^a \left[\kappa_\calS(f(x_0)-f_{\rm low})\epsilon^{-\sfrac{p+1}{p}} +|\calH|\right]
  + \kappa^b + 1
\eeqn
iterations and evaluations of $f$ and its first $p$ derivatives to return a
point $x_\epsilon\in \calF$ such that $\chi_f(x_\epsilon,\epsilon) \leq
\epsilon$.
}

\proof{ If termination occurs at iteration 0, the theorem obviously holds.
  Assume therefore that termination occurs at iteration $k+1$, in which case
  there must be at least one successful iteration.  We may therefore deduce
  the desired bound from Theorem~\ref{a-comp1}, Lemma~\ref{a-succ-unsucc} and the
  fact that each successful iteration involves the evaluation of $f(x_k)$ and
  $\{\nabla_x^if_{\calW_k}(x_k)\}_{i=1}^p$, while each unsuccessful iteration only
  involves that of $f(x_k)$ and $\nabla_x^1f_{\calW_k}(x_k)$.
} 

\noindent
Note that we may count derivatives' evaluations in Theorem~\ref{final-comp}
because only the derivatives of $f_{\calW_k}$ are ever evaluated, and these are
well-defined. For completeness, we state the complexity bound of the important
purely Lipschitzian case.

\lcor{final-comp-Lip}{
Suppose that AS.1-AS.4 hold and $\calH = \emptyset$. Then
Algorithm~\ref{psarp-algo} requires at most
\[
\kappa^a \left[\kappa_\calS(f(x_0)-f_{\rm low})\epsilon^{-\sfrac{p+1}{p}}\right]
  + \kappa^b + 1
\]
iterations and evaluations of $f$ and its first $p$ derivatives to return a
point $x_\epsilon\in \calF$ such that
\[
\chi_f(x_\epsilon)
\eqdef
  \left|\min_{\mystack{x+d\in\calF}{\|d\|\leq 1} }
      \nabla_x^1 f_{\calW(x)}(x)^T d\right|
      \leq \epsilon.
\]
}

\proof{Directly follows from Theorem~\ref{final-comp}, $\calH=\emptyset$ and
  the obesrvation that $\calR(x,\epsilon) = \Re^n$ for all $x\in \calF$ since
  $\calC(x,\epsilon) = \emptyset$.}

\numsection{Evaluation complexity for general convex $\calF$}\label{true-model-s}

The two-sided model \req{mod-hir} has clear advantages, the main ones being
that, except at the origin where it is non-smooth, it is polynomial and has
finite gradients (and higher derivatives) over each of its two branches.  It
is not however without drawbacks.  The first of these is that its prediction
for the gradient (and higher derivatives) is arbitrarily inaccurate as the
origin is approached, the second being its evaluation cost which is typically
higher than evaluating $|x+s|^q$ or its derivative directly.  In particular,
it is the first drawback that required the careful analysis of
Lemma~\ref{a-grad-lemma}, in turn leading, via Lemma~\ref{a-big-grad-l}, to
the crucial Lemma~\ref{a-away-l}. This is significant because this last lemma,
in addition to the use of \req{mod-hir} and the requirement that $p$ must be
odd, also requires the 'kernel-centered' assumption \req{Fbox}, a
sometimes undesirable restriction of the feasible domain geometry.

In the case where evaluating $f_\calN$ is very expensive and the convex
$\calF$ is not 'kernel-centered', it may sometimes be
acceptable to push the difficulty of handling the non-Lipschitzian nature of
the $\ell_q$ norm regularization in the subproblem of computing $s_k$, if evaluations of
$f_\calN$ can be saved.  In this context, a simple alternative is then to use
\beqn{mod-hir-true}
m_i(x_i,s_i) = |x_i+s_i|^q
\tim{for}
i \in \calH
\eeqn
that is $m_i(x_i,s_i) = f_i(x_i+s_i)$ for $i \in \calH$. The cost of finding a
suitable step satisfying \req{mod-term} may of course be increased, but, as we already
noted, this cost is irrelevant for worst-case evaluation analysis as long as
only the evaluation of $f_\calN$ and its derivatives is taken into account.
The choice \req{mod-hir-true} clearly maintains the overestimation property of
Lemma~\ref{a-overestimation-l}.  Moreover, it is easy to verify (using AS.3
and \req{mod-hir-true}) that
\beqn{good-grad}
\|\nabla_x f_{\calW_k^+}(x_k+s_k) - \nabla_s^1m_{\calW_k^+}(x_k,s_k)\|
= \|\nabla_x f_\calN(x_k+s_k) - \nabla_s^1m_\calN(x_k,s_k)\|
\leq L_{\max} \|s_k\|^p.
\eeqn
This in turn implies that the proof of Lemma~\ref{a-schi+-lemma} can be extended
without requiring \req{heartguy} and using $\calO_* =
\calH\setminus\calC_k^+$.  The derivation of \req{a-errdqm} then simplifies
because of \req{good-grad} and holds for all $i\in\calH\setminus\calC_k^+$
with $L(\alpha) = L_{\max}$, so that \req{a-schi+} holds for all $k\in\calS$,
the assumption \req{eps-small-1} being now irrelevant.  This result then
implies that the distinction made between $\calS_\heartsuit$,
$\calS_\diamondsuit$, $S\cal_\clubsuit$ and  $\calS_{\|s\|}$ is
unecessary because \req{a-schi+} holds for all $k\in\calS= \calS_\heartsuit$.
Moreover, since we no longer need Lemma~\ref{a-away-l} to prove that
$\calS_\clubsuit = \emptyset$, we no longer need the restrictions that $p$ is
odd and \req{Fbox} either.  As consequence, we deduce that
Theorem~\ref{a-comp1} holds for arbitrary $p\geq 1$ and for arbitrary convex,
closed non-empty $\calF$, without the need to assume \req{a-eps-upper} and
with $L(\alpha)$ replaced by $L_{\max}$ in \req{a-kapS-def}. Without altering
Lemma~\ref{a-succ-unsucc}, we may therefore deduce the following complexity
result.

\lthm{final-comp-true}{
Suppose that AS.1, AS.2 (without the restriction that $p$ must be odd), AS.3
and AS.4 hold. Then
Algorithm~\ref{psarp-algo} using the true models \req{mod-hir-true} for
$i\in\calH$ requires at most
\[
\kappa^a \left[\kappa_\calS^{\rm true}(f(x_0)-f_{\rm low})\epsilon^{-\sfrac{p+1}{p}} +|\calH|\right]
  + \kappa^b + 1
\]
iterations and evaluations of $f_\calN$ and its first $p$ derivatives to return a
point $x_\epsilon\in \calF$ such that $\chi_f(x_\epsilon,\epsilon) \leq
\epsilon$, where
\[
\kappa_\calS^{\rm true} = \bigfrac
      {(p+1)!}{\eta \, \sigma_{\min}\, \varsigma_{\min}^{p+1} }
 \Big[2 |\calN|\varsigma^{p+1}\left(L+\theta+\bigfrac{\gamma_2}{p!}\right) \Big]^{\frac{p+1}{p}}.
\]
}

\noindent
As indicated, the complexity is expressed in this theorem in terms of
evaluations of $f_\calN$ and its derivatives only.  The evaluation count for
the terms $f_i$ ($i\in \calH$) may be higher since these terms are
evaluated in computing the step $s_k$ using the models \req{mod-hir-true}.
Note that the difficulty of handling infinite derivatives is passed on to the
subproblem solver in this approach.

Moreover, it also results from the analysis in this section that one may
consider objective functions of the form
\[
f(x) = f_\calN(x) + f_\calH(x)
\]
and prove an $O(\epsilon^{-\frac{p+1}{p}})$ evaluation compexity bound if
$f_\calN$ has Lipschitz continuous derivatives of order $p$ and if
$m_\calH(x_k,s) = f_\calH(x_k+s)$, passing all difficulties associated with
$f_\calH$ to the subproblem of computing the step $s_k$.

As it turns out, an evaluation complexity bound may also be computed if one
insist on using the Taylor's models \req{mod-hir} while allowing the feasible
set to be an arbitrary convex, closed and non-empty set.  Not surprisingly,
the bound is (significantly) worse than that provided by
Theorem~\ref{final-comp}, but has the merit of existing.  Its derivation is
based on the observation that \req{a-Tayl-rest} in
Lemma~\ref{a-overestimation-l} and \req{Pik-bound} imply that, for $i
\in\calH\setminus\calC_k^+$,
\beqn{very-bad-grad-err}
|\nabla_{s_i}^1|x_i+s_i|^q -\nabla_{s_i}^1m_i(x_i,s_i)|
\leq q \left(\min\Big[|x_i|,|x_i+s_i|\Big]\right)^{q-p-1} |\mu(x_i,s_i)|^p
\leq q \epsilon^{q-p-1} |s_i|^p.
\eeqn
This bound can then be used in a variant of Lemma~\ref{a-schi+-lemma} just
like \req{good-grad} was in Section~\ref{true-model-s}.  In the updated
version of Lemma~\ref{a-schi+-lemma}, we replace $L(\alpha)$ by
\[
L_* \eqdef |\calN|\,\varsigma_{\max}^p\,L_{\max} +
|\calH|\,\varsigma_{\max}^p\,q
\]
and \req{a-errdqm} now becomes
\[
\|\nabla_x^1f_{\calW_k^+}(x_{k+1}) -\nabla_s^1m_{\calW_k^+}(x_k,s_k)\|
\leq \left[L_* \epsilon^{q-p-1}+\bigfrac{|\calN|}{p!}\varsigma^p\sigma_{\max}\right]\|s_k\|^p.
\]
This results in replacing \req{a-con2-ARC2CC-sl-1b} by
\beqn{ns-2}
\chi_f(x_{k+1})
\leq (L_*\epsilon^{q-p-1}+\theta + \bigfrac{|\calN|}{p!}\,
     \varsigma_{\max}^{p+1}\,\sigma_{\max})\|s_k\|^p
\leq (L_*+\theta + \bigfrac{|\calN|}{p!}\,
     \varsigma_{\max}^{p+1}\,\sigma_{\max})\epsilon^{q-p-1}\|s_k\|^p
\eeqn
and therefore \req{a-schi+} is replaced by
\[
\|s_k\|
\geq \left[2\left(L_* +\theta+\bigfrac{|\calN|}{p!}\,
         \varsigma_{\max}^{p+1}\,\sigma_{\max}\right)\right]^{-\frac{1}{p}} \epsilon^{\frac{p+2-q}{p}}.
\]
We may now follow the steps leading to Theorem~\ref{final-comp-true} and deduce a new
complexity bound.

\lthm{final-comp-worse}{
Suppose that AS.1--AS.4 hold. Then
Algorithm~\ref{psarp-algo} using the Taylor models \req{mod-hir} for
$i\in\calH$ requires at most
\[
\kappa^a \left[\kappa_\calS^*(f(x_0)-f_{\rm low})\epsilon^{-\sfrac{(p+2-q)(p+1)}{p}} +|\calH|\right]
  + \kappa^b + 1
\]
iterations and evaluations of $f$ and its first $p$ derivatives to return a
point $x_\epsilon\in \calF$ such that $\chi_f(x_\epsilon,\epsilon) \leq
\epsilon$, where
\[
\kappa_\calS^* = \bigfrac
      {(p+1)!}{\eta \, \sigma_{\min}\, \varsigma_{\min}^{p+1} }
      \Big[2\left( L_*+\theta + \bigfrac{|\calN|}{p!}\,
        \varsigma_{\max}^{p+1}\,\gamma_2\right) \Big]^{\frac{p+1}{p}}.
\]
}

\noindent
Observe that, due to the second inequality of \req{ns-2}, $\theta$ can be
replaced in \req{mod-term} by $\theta_* = \theta \epsilon^{q-p-1}$, making the
termination condition for the step computation very weak.

\numsection{Further discussion}\label{discuss-s}

The above results suggest some additional comments.
\begin{itemize}
\item The complexity result in $O(\epsilon^{-(p+1)/p})$ evaluations obtained
  in Theorem~\ref{final-comp} is identical in order to that presented in
  \cite{BirgGardMartSantToin17} for the unstructured unconstrained
  and in \cite{CartGoulToin16a} for the unstructured convexly constrained
  cases. It is remarkable that incorporating non-Lipschitzian singularities
  in the objective function does not affect the worst-case evaluation
  complexity of finding an $\epsilon$-approximate first-order critical point.

\item Interestingly, Corollary~\ref{final-comp-Lip} also shows that using
  partially separable structure does not affect the evaluation complexity
  either, therefore allowing cost-effective use of problem structure with
  high-order models.

\item The algorithm\footnote{And theory, if one restricts one's attention to
  the case where $\calH  = \emptyset$.} presented here is considerably simpler
  than that discussed in \cite{ConnGoulSartToin96a,ConnGoulToin00} in the
  context of structured trust-regions. In addition, the present assumptions
  are also weaker.  Indeed, an additional condition on long steps
  (see AA.1s in \cite[p.364]{ConnGoulToin00}) is no longer needed.

\item Can one use even order models with Taylor models in the present
  framework? The main issue is that, when $p$ is even, the two-sided model
  $T_{|\cdot|^q,p}(x_i,s_i)$ is no longer always an overestimate of
  $|x_i+s_i|^q$ when $|x_i+s_i|>|x_i|$, as can be verified from
  \req{a-Tayl-rest}. While this can be taken care of by adding a
  regularization term to $m_i$, the necessary size of the regularization
  parameter may be unbounded when the iterates are sufficiently close from the
  singularity.  This in turn destroys the good complexity
  because it forces the algorithm to take much too short steps.

  An alternative is to use mixed-orders models, that is models of even order
  ($p$, say) for the $f_i$ whose index is in $\calN$ and odd order models for those
  with index in $\calH$. However, this last (odd) order has to be at least as
  large as $p$, because it is the lowest order which dominates in the crucial
  Lemma~\ref{a-schi+-lemma} where the length is bounded below away from the
  singularity. The choice of a $(p+1)$-st order model for $i\in \calH$ is then
  most natural.

\item A variant of the algorithm can be stated where it is
  possible for a particular $x_i$ to leave the $\epsilon$-neighbourhood of
  zero, provided the associated step results in a significant (in view of
  Theorem~\ref{a-comp1}) objective function decrease, such as a multiple of
  $\epsilon^{(p+1)/p}$ or some $\epsilon$-independent constant.  These
  decreases can then be counted separately in the argument of Theorem~\ref{a-comp1}
  and cycling is impossible since there can be only a finite number of such
  decreases.
\end{itemize}

\numsection{Conclusions}\label{concl-s}

We have considered the problem of minimizing a partially-separable nonconvex
objective function $f$ involving non-Lipschitzian $q$-norm regularization
terms and subject to general convex constraints. Problems of this type are
important in many areas, including data compression, image processing and
bioinformatics. We have shown that the introduction of the non-Lipschitzian
singularities and the exploitation of problem structure do not affect the
worst-case evaluation complexity. More precisely, we have first defined
$\epsilon$-approximate first-order critical points for the considered class of
problems in a way that make the obtained complexity bounds comparable to
existing results for the purely Lipschitzian case. We have then shown that, if
$p$ is the (odd) degree of the models used by the algorithm, if the feasible
set is 'kernel-centered' and if Taylor models are used for the $q$-norm
regularization terms, the bound of $O(\epsilon^{-\frac{p+1}{p}})$ evaluations
of $f$ and its relevant derivatives (derived for the Lipschitzian case in
\cite{CartGoulToin16a}) is preserved in the presence of non-Lipschitzian
singularities.  In addition, we have shown that partially-separable structure
present in the problem can be exploited (especially for high degree derivative
tensors) without affecting the evaluation complexity either. We have also
shown that, if the difficulty of handling the non-Lipschitzian regularization
terms is passed to the subproblem (which can be meaningfull if evaluating the
other parts of the objective function is very expensive) in that non-Lipschitz
models are used for these terms, then the same bounds hold in terms of
evaluation of the expensive part of the objective function, without the
restriction that the feasible set be `kernel-centered'. A worse complexity
bound has finally been provided in the case where one uses Taylor models for the
$q$-norm regularization terms with a general convex feasible set.

These objectives have been attained by introducing a new first-order
criticality measure as well as the new two-sided model of the singularity
given by \req{modiH}, which exploits the inherent symmetry and provides a
useful overestimate of the $|x|^q$ if its order is chosen odd, without the
need for smoothing functions.


\section*{\footnotesize Acknowledgements}

{\footnotesize

Xiaojun Chen would like to thank Hong Kong Research Grant Council for grant
PolyU153000/15p. Philippe Toint would like to thank the Belgian Fund for Scientific
Research (FNRS), the University of Namur and the Hong Kong Polytechnic
University for their support while this research was being conducted.

\bibliographystyle{plain}
\bibliography{/home/pht/bibs/refs}
}
\end{document}